\documentclass[12pt]{article}
\usepackage{amsmath}
\usepackage{amssymb}
\renewcommand\smallskip{\vskip\smallskipamount}
\renewcommand\medskip{\vskip\medskipamount}
\renewcommand\bigskip{\vskip\bigskipamount}

\begin{document}

\begin{center}
\begin{large}
\textbf{The Local Isometric Embedding in} $\mathbb{R}^{3}$
\textbf{of Two-Dimensional Riemannian Manifolds With Gaussian
Curvature Changing Sign to Finite Order on a Curve}
\end{large}

\bigskip\medskip
MARCUS A. KHURI \footnote{The author was partially supported by an
NSF Postdoctoral Fellowship and NSF Grant DMS-0203941.}

\bigskip
\bigskip
\bigskip
\begin{abstract}
We consider two natural problems arising in geometry which are
equivalent to the local solvability of specific equations of
Monge-Amp\`{e}re type.  These two problems are: the local
isometric embedding problem for two-dimensional Riemannian
manifolds, and the problem of locally prescribed Gaussian
curvature for surfaces in $\mathbb{R}^{3}$.  We prove a general
local existence result for a large class of Monge-Amp\`{e}re
equations in the plane, and obtain as corollaries the existence of
regular solutions to both problems, in the case that the Gaussian
curvature vanishes to arbitrary finite order on a single smooth
curve.
\end{abstract}
\bigskip
\bigskip

\textbf{0.  Introduction}
\end{center}\smallskip

   Let $(M^{2},ds^{2})$ be a two-dimensional Riemannian manifold.
A well-known problem is to ask when can one realize this, locally,
as a small piece of a surface in $\mathbb{R}^{3}$.  This question
has only been partially answered.\par
   Suppose that the first fundamental form,
$ds^{2}=Edu^{2}+2Fdudv+Gdv^{2}$, is given in the neighborhood of a
point, say $(u,v)=0$.  Let $K$ be the Gaussian curvature, then the
known results are as follows.  The question is answered
affirmatively in the case that $ds^{2}$ is analytic or $K(0)\neq
0$; these classical results can be found in [8], [16], and [17].
In the case that $K\geq 0$ and $ds^{2}$ is sufficiently smooth, or
$K(0)=0$ and $\nabla K(0)\neq 0$, C.-S. Lin provides an
affirmative answer in [12] and [13] (a simplified proof of the
later result has been given by Q. Han [4]).  If $K\leq 0$ and
$\nabla K$ possesses a certain nondegeneracy, Han, Hong, and Lin
[6] show that an embedding always exists.  Furthermore, if
$(u,v)=0$ is a nondegenerate critical point for $K$ and $ds^{2}$
is sufficiently smooth, then the author provides an affirmative
answer in [11]. However, A. V. Pogorelov has given a
counterexample in [15], where he constructs a $C^{2,1}$ metric
with no $C^{2}$ isometric embedding in $\mathbb{R}^{3}$. More
recently, other counterexamples for metrics with low regularity
have been proposed by Nadirashvili and Yuan [14], and local
nonexistence results for smooth Monge-Amp\`{e}re equations have
been obtained in [10].  In this paper we prove the
following,\medskip

\textbf{Theorem 0.1.}  \textit{Let $ds^{2}\in C^{r}$, $r\geq 60$,
and suppose that $\sigma$ is a geodesic passing through the
origin.  If $K$ vanishes to finite order on $\sigma$, then there
exists a $C^{r-36}$ local isometric embedding into
$\mathbb{R}^{3}$.}\medskip

\textbf{Remark.}  \textit{The geodesic hypothesis on $\sigma$ is
actually unnecessary, and is only included so that Theorem} 0.1
\textit{arises as a corollary of our main result, Theorem} 0.3
\textit{below.  Please see the appendix for the justification.
Also, a similar result has been obtained independently by Q. Han}
[5].\medskip

   We begin by deriving the appropriate equations for study.  Our
goal is to find three functions $x(u,v)$, $y(u,v)$, $z(u,v)$, such
that $ds^{2}=dx^{2}+dy^{2}+dz^{2}$.  The following strategy was
first used by J. Weingarten [21].  We search for a function
$z(u,v)$, with $|\nabla z|$ sufficiently small, such that
$ds^{2}-dz^{2}$ is flat in a neighborhood of the origin.  Suppose
that such a function exists, then since any Riemannian manifold of
zero curvature is locally isometric to Euclidean space (via the
exponential map), there exists a smooth change of coordinates
$x(u,v)$, $y(u,v)$ such that $dx^{2}+dy^{2}=ds^{2}-dz^{2}$, that
is, $ds^{2}=dx^{2}+dy^{2}+dz^{2}$.  Therefore, our problem is
reduced to finding $z(u,v)$ such that $ds^{2}-dz^{2}$ is flat in a
neighborhood of the origin.  A computation shows that this is
equivalent to the local solvability of the following equation,
\begin{eqnarray}
& &(z_{11}-\Gamma^{i}_{11}z_{i})(z_{22}-\Gamma^{i}_{22}z_{i})
-(z_{12}-\Gamma^{i}_{12}z_{i})^{2}\\
&=& K(EG-F^{2}-Ez^{2}_{2}-Gz^{2}_{1} +2Fz_{1}z_{2}),\nonumber
\end{eqnarray}
where $z_{1}=\partial z/\partial u$, $z_{2}=\partial z/\partial
v$, $z_{ij}$ are second derivatives of $z$, and $\Gamma_{jk}^{i}$
are Christoffel symbols.\par
   Equation (1) is a second order Monge-Amp\`{e}re equation.
Another well-known and related problem, which is equivalent to the
local solvability of a second order Monge-Amp\`{e}re equation, is
that of locally prescribing the Gaussian curvature for surfaces in
$\mathbb{R}^{3}$.  That is, given a function $K(u,v)$ defined in a
neighborhood of the origin, when does there exist a piece of a
surface $z=z(u,v)$ in $\mathbb{R}^{3}$ having Gaussian curvature
$K$?  This problem is equivalent to the local solvability of the
equation
\begin{equation}
z_{11}z_{22}-z_{12}^{2}=K(1+|\nabla z|^{2})^{2}.
\end{equation}
For this problem we obtain a similar result to that of Theorem
0.1.\medskip

\textbf{Theorem 0.2.}  \textit{Let $\sigma$ be a smooth curve
passing through the origin.  If $K\in C^{r}$, $r\geq 58$, and $K$
vanishes to finite order on $\sigma$, then there exists a piece of
a $C^{r-34}$ surface in $\mathbb{R}^{3}$ with Gaussian curvature
$K$.}\medskip

   With the goal of treating both problems simultaneously, we will
study the local solvability of the following general
Monge-Amp\`{e}re equation
\begin{equation}
\det(z_{ij}+a_{ij}(u,v,z,\nabla z))=Kf(u,v,z,\nabla z),
\end{equation}
where $a_{ij}(u,v,p,q)$ and $f(u,v,p,q)$ are smooth functions of
$p$ and $q$, $f>0$, $K$ vanishes to finite order along a smooth
curve $\sigma$ passing through the origin, and $a_{ij}$ vanishes
along $\sigma$ to an order greater than or equal to one degree
less than that of $K$.  Clearly equation (2) is of the form (3),
and equation (1) is of the form (3) if $\Gamma_{jk}^{i}$ vanishes
to the order of one degree less than that of $K$ along $\sigma$,
which we assume without loss of generality. More precisely, since
$\sigma$ is a geodesic we can introduce geodesic parallel
coordinates, such that $\sigma$ becomes the $v$-axis and
$ds^{2}=du^{2}+h^{2}dv^{2}$, for some $h\in C^{r-1}$ satisfying
\begin{equation*}
h_{uu}=-Kh,\text{ }\text{ }h(0,v)=1,\text{ }\text{ }h_{u}(0,v)=0.
\end{equation*}
It then follows that the Christoffel symbols vanish to the
appropriate order along the $v$-axis.  We will prove\medskip

\textbf{Theorem 0.3.}  \textit{Let $\sigma$ be a smooth curve
passing through the origin.  If $K$, $a_{ij}$, $f\in C^{r}$,
$r\geq 58$, $K$ vanishes to finite order along $\sigma$, and
$a_{ij}$ vanishes to an order greater than or equal to one degree
less than that of $K$ along $\sigma$, then there exists a
$C^{r-34}$ local solution of} (3).\medskip

   Equation (3) is elliptic if $K>0$, hyperbolic if $K<0$, and
of mixed type if $K$ changes sign in a neighborhood of the origin.
If $K(0)=0$ and $\nabla K(0)\neq 0$ [13], then (3) is a nonlinear
type of the Tricomi equation.  While if the origin is a
nondegenerate critical point for $K$ [11], then (3) is a nonlinear
type of Gallerstedt's equation [3].  In our case, assuming that
$K$ vanishes to some finite order $n+1\in\mathbb{Z}_{>0}$ along
$\sigma$ (ie. all derivatives up to and including order $n$ vanish
along $\sigma$), and $a_{ij}$ vanishes at least to order $n$ along
$\sigma$, the linearized equation for (3) may be put into the
following canonical form after adding suitable first and second
order perturbation terms and making an appropriate change of
coordinates,
\begin{equation}
Lu=y^{n+1}A_{1}u_{xx}+u_{yy}+y^{n-1}A_{2}u_{x}+A_{3}u_{y}+A_{4}u,
\end{equation}
where the $A_{i}$ are smooth functions and $A_{1}>0$ or $A_{1}<0$.
It will be shown that this special canonical form is amenable to
the making of estimates, even in the case that (4) changes type
along the line $y=0$.\par
   From now on we assume that $n>0$ is even, since the case when $n$
is odd may be treated by the results in [12] and [6] where $K$ is
assumed to be nonnegative or nonpositive, and the case $n=0$ may
be treated by the methods of [13].  Furthermore, we assume without
loss of generality that the curve $\sigma$ is given by an equation
$\widetilde{H}(u,v)=0$, where $\widetilde{H}\in C^{\infty}$ and
$\widetilde{H}_{v}|_{\sigma}\geq M_{1}$ for some constant
$M_{1}>0$.  Let $\varepsilon$ be a small parameter and set
$u=\varepsilon^{2}x$, $v=\varepsilon^{2}y$,
$z=u^{2}/2+\varepsilon^{5}w$ (the $x,y$ used here are not the same
as those appearing in (4)). Substituting into (3), we obtain
\begin{equation}
\Phi(w):=(1+\varepsilon w_{xx}+a_{11})(\varepsilon
w_{yy}+a_{22})-(\varepsilon w_{xy}+a_{12})^{2}-Kf=0.
\end{equation}
By the assumptions of Theorem 0.3 we may write
\begin{equation*}
a_{ij}=\varepsilon^{2n}H^{n}(x,y)P_{ij}(\varepsilon, x,y,w,\nabla
w)
\end{equation*}
and
\begin{equation*}
Kf=\varepsilon^{2(n+1)}H^{n+1}(x,y)P(\varepsilon, x,y,w,\nabla w),
\end{equation*}
where
\begin{equation*}
H=\varepsilon^{-2}\widetilde{H},\text{ }\text{ }\text{ }\text{
}\text{ }H_{y}|_{\sigma}\geq M_{1},\text{ }\text{ }\text{ }\text{
}\text{ }P\geq M_{2}
\end{equation*}
for some constant $M_{2}>0$ independent of $\varepsilon$, and
$P_{ij}$, $P$ are $C^{r}$ with respect to $x,y$ and $C^{\infty}$
with respect to the remaining variables.  Then (5) becomes
\begin{eqnarray}
\Phi(w)&=&(1+\varepsilon
w_{xx}+\varepsilon^{2n}H^{n}P_{11})(\varepsilon
w_{yy}+\varepsilon^{2n}H^{n}P_{22})\nonumber\\
& &-(\varepsilon
w_{xy}+\varepsilon^{2n}H^{n}P_{12})^{2}-\varepsilon^{2(n+1)}H^{n+1}P\\
&=&0.\nonumber
\end{eqnarray}
Choose $x_{0},y_{0}>0$ and define the rectangle $X=\{(x,y)\mid
|x|<x_{0},|y|<y_{0}\}$.  Then solving $\Phi(w)=0$ in $X$, is
equivalent to solving (3) locally at the origin.\par
   In the following sections, we shall study the linearization of
(6) about some function $w$.  In section $\S 1$ the linearization
will be reduced to the canonical form (4). Existence and
regularity for the modified linearized equation will be obtained
in section $\S 2$.  In section $\S 3$ we make the appropriate
estimates in preparation for the Nash-Moser iteration procedure.
Finally, in $\S 4$ we apply a modified version of the Nash-Moser
procedure and obtain a solution of (6).  An appendix is included
in section $\S 5$ in order to justify removing the geodesic
hypothesis from Theorem 0.1.
\bigskip
\begin{center}
\textbf{1.  Reduction to Canonical Form}
\end{center}
\bigskip

   In this section we will bring the linearization of (6) into
the canonical form (4).  This shall be accomplished by adding
certain perturbation terms and making appropriate changes of
variables.  The process will entail defining a sequence of linear
operators $L_{i}$, $1\leq i\leq 7$, where $L_{1}$ is the
linearization of (6) and $L_{7}$ is of the form (4); furthermore,
$L_{i+1}$ will differ from $L_{i}$ by a perturbation term or by a
change of variables.\par
   Fix a constant $C>0$, and let $w\in C^{\infty}(\mathbb{R}^{2})$
be such that $|w|_{C^{16}}\leq C$.  Then the linearization of (6)
evaluated at $w$ is given by
\begin{equation}
L_{1}(w)=\sum_{i,j}b_{ij}^{1}\partial_{x_{i}x_{j}}+\sum_{i}b_{i}^{1}
\partial_{x_{i}}+b^{1},
\end{equation}
where $x_{1}=x$, $x_{2}=y$ and
\begin{eqnarray*}
b_{11}^{1}&=&\varepsilon(\varepsilon w_{yy}+
\varepsilon^{2n}H^{n}(x,y)P_{22}(\varepsilon,x,y,w,\nabla w)),\\
b_{12}^{1}=b_{21}^{1}&=&-\varepsilon(\varepsilon w_{xy}+
\varepsilon^{2n}H^{n}(x,y)P_{12}(\varepsilon,x,y,w,\nabla w)),\\
b_{22}^{1}&=&\varepsilon(1+\varepsilon w_{xx}+
\varepsilon^{2n}H^{n}(x,y)P_{11}(\varepsilon,x,y,w,\nabla w)),\\
b_{1}^{1}&=&\varepsilon^{2n}H^{n}(x,y)P_{1}(\varepsilon,x,y,w,\nabla w),\\
b_{2}^{1}&=&\varepsilon^{2n}H^{n}(x,y)P_{2}(\varepsilon,x,y,w,\nabla w),\\
b^{1}&=&\varepsilon^{2n}H^{n}(x,y)P_{3}(\varepsilon,x,y,w,\nabla
w),
\end{eqnarray*}
for some $P_{1}$, $P_{2}$, $P_{3}$.  If $\varepsilon$ is
sufficiently small, we may solve for $\varepsilon
w_{yy}+\varepsilon^{2n}H^{n}P_{22}$ in equation (6) to obtain
\begin{equation}
\varepsilon w_{yy}+\varepsilon^{2n}H^{n}P_{22}=\frac{1}
{1+\varepsilon Q}[(\varepsilon
w_{xy}+\varepsilon^{2n}H^{n}P_{12})^{2}+\varepsilon^{2(n+1)}H^{n+1}P
+\Phi(w)],
\end{equation}
where $Q(\varepsilon,x,y,w,\nabla
w,\nabla^{2}w)=w_{xx}+\varepsilon^{2n-1}H^{n}P_{11}$.  Plugging
(8) into (7) we have,
\begin{eqnarray*}
L_{2}(w)&:=&L_{1}(w)-\frac{\varepsilon\Phi(w)}
{1+\varepsilon Q}\partial_{xx}\\
&=&\sum_{i,j}b_{ij}^{2}\partial_{x_{i}x_{j}}+\sum_{i}b_{i}^{2}
\partial_{x_{i}}+b^{2},
\end{eqnarray*}
where
\begin{equation*}
b_{11}^{2}=\frac{\varepsilon(\varepsilon
w_{xy}+\varepsilon^{2n}H^{n}P_{12})^{2}+\varepsilon^{2n+3}H^{n+1}P}
{1+\varepsilon Q}.
\end{equation*}
Next define $L_{3}(w)$ by,
\begin{eqnarray}
L_{3}(w)&:=&\frac{1}{\varepsilon(1+\varepsilon
Q)}L_{2}(w)\\
&=&\sum_{i,j}b_{ij}^{3}\partial_{x_{i}x_{j}}+\sum_{i}b_{i}^{3}
\partial_{x_{i}}+b^{3}.\nonumber
\end{eqnarray}\par
   To simplify (9), we will make a change of variables that will
eliminate the mixed second derivative term.  In constructing this
change of variables we will make use of the following lemma from
ordinary differential equations.\medskip

\textbf{Lemma 1.1 [1].}  \textit{Let $G(x,t)$ be a $C^{l}$ real
valued function in the closed rectangle $|x-s|\leq T_{1}$,
$|t|\leq T_{2}$.  Let $T=\sup|G(x,t)|$ in this domain.  Then the
initial value problem $dx/dt=G(x,t)$, $x(0)=s$, has a unique
$C^{l+1}$ solution defined on the interval
$|t|\leq\min(T_{2},T_{1}/T)$.  Moreover, $x(s,t)$ is $C^{l}$ with
respect to $s$.}\medskip

   We now construct the desired change of variables.  For any
domain $\Omega\subset\mathbb{R}^{2}$, and constant $\mu$, let $\mu
\Omega=\{\mu(x,y)\mid(x,y)\in \Omega\}$.\medskip

\textbf{Lemma 1.2.}  \textit{For $\varepsilon$ sufficiently small,
there exists a $C^{r}$ diffeomorphism}
\begin{equation*}
\xi=\xi(x,y),\text{   }\eta=y,
\end{equation*}
\textit{of a domain $X_{1}$ onto $\mu_{1}X$, where $\mu_{1}>1$,
such that in the new variables $(\xi,\eta)$, $L_{3}(w)$ is denoted
by $L_{4}(w)$ and is given by}
\begin{equation*}
L_{4}(w)=\sum_{i,j}b_{ij}^{4}\partial_{x_{i}x_{j}}+\sum_{i}b_{i}^{4}
\partial_{x_{i}}+b^{4},
\end{equation*}
\textit{where $x_{1}=\xi$, $x_{2}=\eta$, and}
\begin{eqnarray*}
b_{11}^{4}\!\!\!\!\!&=&\!\!\!\!\varepsilon^{2(n+1)}H^{n+1}P_{11}^{4},\\
b_{12}^{4}=b_{21}^{4}\!\!\!\!\!\!&\equiv&\!\!\! 0,\\
b_{22}^{4}\!\!\!\!\!&\equiv&\!\!\!\! 1,\\
b_{1}^{4}\!\!\!\!\!&=&\!\!\!\!\varepsilon^{2n}H^{n}P_{1}^{41}\!+n\varepsilon^{2n}H^{n-1}P_{1}^{42}\!
+[\partial_{x}(\frac{\Phi(w)}{2(1\!+\varepsilon
Q)^{2}})+\frac{\partial_{x}\Phi(w)} {2(1\!+\varepsilon
Q)^{2}}]\xi_{x},\\
b_{2}^{4}\!\!\!\!\!&=&\!\!\!\! b_{2}^{3},\\
b^{4}\!\!\!\!\!&=&\!\!\!\! b^{3},
\end{eqnarray*}
\textit{for some $P_{11}^{4}$, $P_{1}^{41}$, $P_{1}^{42}$, and
$P_{11}^{4}\geq C_{1}$ for some constant $C_{1}>0$ independent of
$\varepsilon$ and $w$.  Furthermore
$\sum|b_{ij}^{4}|_{C^{12}}+|b_{i}^{4}|_{C^{12}}+|b^{4}|_{C^{12}}\leq
C_{2}$, for some $C_{2}$ independent of $\varepsilon$ and
$w$.}\medskip

\textit{Proof.}  Using the chain rule we find that
$b_{12}^{4}=b_{12}^{3}\xi_{x}+b_{22}^{3}\xi_{y}$. Therefore, we
seek a smooth function $\xi(x,y)$ such that
\begin{equation}
b_{12}^{4}=b_{12}^{3}\xi_{x}+b_{22}^{3}\xi_{y}=0\text{ }\text{
}\text{ in }\text{ }X_{1},\text{ }\text{ }\text{ }\xi(x,0)=x,
\end{equation}
where $X_{1}$ will be defined below.  Since $b_{22}^{3}\equiv 1$,
the line $y=0$ will be non-characteristic for (10).  Then by the
theory of first order partial differential equations, (10) is
reduced to the following system of first order ODE:
\begin{eqnarray*}
\dot{x}&=&b_{12}^{3},\text{ }x(0)=s,\text{ }-\mu_{1}x_{0} \leq
s\leq\mu_{1}x_{0},\\
\dot{y}&=&1,\text{ }\text{ }\text{ }y(0)=0,\\
\dot{\xi}&=&0,\text{ }\text{ }\text{ }\xi(0)=s,
\end{eqnarray*}
where $x=x(t)$, $y=y(t)$, $\xi(t)=\xi(x(t),y(t))$ and $\dot{x}$,
$\dot{y}$, $\dot{\xi}$ are derivatives with respect to $t$.\par
   Choose $\mu_{1}>1$.  We first show that the characteristic curves, given
parametrically by $(x,y)=(x(t),t)$, exist globally for
$-\mu_{1}y_{0}\leq t\leq\mu_{1}y_{0}$.  We apply Lemma 1.1 with
$T_{1}=2\mu_{1}x_{0}$, and $T_{2}=\mu_{1}y_{0}$, to the
initial-value problem $\dot{x}=b_{12}^{3}$, $x(0)=s$.  Let $T$ be
as in Lemma 1.1.  Since $|w|_{C^{16}}\leq C$, we have
\begin{equation*}
T=\sup_{X_{1}}|b_{12}^{3}|\leq\varepsilon C_{3},
\end{equation*}
for some $C_{3}$ independent of $\varepsilon$.  Then for
$\varepsilon$ small, $T\leq\frac{2x_{0}}{y_{0}}$, implying that
\begin{equation*}
\min(T_{2},T_{1}/T)=\mu_{1}y_{0}.
\end{equation*}
Then Lemma 1.1 gives the desired global existence.\par
   Let $X_{1}$ be the domain with boundary
consisting of the two lines $y=\pm\mu_{1}y_{0}$, and the two
characteristics passing through $\pm\mu_{1}x_{0}$.  Then the
mapping $(\xi,\eta)$ takes $\partial X_{1}$ onto
$\partial\mu_{1}X$.  We now show that the map
$\rho:\mu_{1}X\rightarrow X_{1}$ given by
$(s,t)\mapsto(x(s,t),y(s,t))=(x(s,t),t)$, is a diffeomorphism.  It
will then follow that the map
$(x,y)\mapsto(\xi(x,y),\eta(x,y))=(s(x,y),y)=\rho^{-1}(x,y)$ is a
diffeomorphism of $X_{1}$ onto $\mu_{1}X$.  To show that $\rho$ is
1-1, suppose that $\rho(s_{1},t_{1})=\rho(s_{2},t_{2})$. Then
$t_{1}=t_{2}$ and $x(s_{1},t_{1})=x(s_{2},t_{2})$, which implies
that $s_{1}=s_{2}$ by uniqueness for the initial-value problem for
ordinary differential equations.  To show that $\rho$ is onto,
take an arbitrary point $(x_{1},y_{1})\in X_{1}$, then we will
show that there exists $s\in[-\mu_{1}x_{0},\mu_{1}x_{0}]$ such
that $\rho(s,y_{1})=(x(s,y_{1}),y_{1})=(x_{1},y_{1})$.  Since the
map
\begin{equation*}
x(s,y_{1}):[-\mu_{1}x_{0},\mu_{1}x_{0}]\rightarrow[x(-\mu_{1}x_{0},y_{1}),
x(\mu_{1}x_{0},y_{1})]
\end{equation*}
is continuous, and $x(-\mu_{1}x_{0},y_{1})\leq x_{1}\leq
x(\mu_{1}x_{0},y_{1})$ by definition of $X_{1}$, the intermediate
value theorem guarantees that there exists
$s\in[-\mu_{1}x_{0},\mu_{1}x_{0}]$ with $x(s,y_{1})=x_{1}$.
Therefore, $\rho$ has a well-defined inverse
$\rho^{-1}:X_{1}\rightarrow\mu_{1}X$.\par
   To show that $\rho^{-1}$ is smooth it is sufficient, by the
inverse function theorem, to show that the Jacobian of $\rho$ does
not vanish at each point of $\mu_{1}X$.  Since
\begin{equation*}
D\rho=\left(%
\begin{array}{cc}
 x_{s} & x_{t}\\
 0 & 1\\
\end{array}%
\right),
\end{equation*}
this is equivalent to showing that $x_{s}$ does not vanish in
$\mu_{1}X$.  Differentiate the equation for $x$ with respect to
$s$ to obtain, $\frac{d}{dt}(x_{s})=(b_{12}^{3})_{x}x_{s}$,
$x_{s}(0)=1$.  Then by the mean value theorem
\begin{equation*}
|x_{s}(s,t)-1|=|x_{s}(s,t)-x_{s}(s,0)|\leq\mu_{1}y_{0}\sup_{X_{1}}
|(b_{12}^{3})_{x}|\sup_{\mu_{1}X}|x_{s}|
\end{equation*}
for all $(s,t)\in\mu_{1}X$.  Thus, since $|w|_{C^{16}}\leq C$,
\begin{equation*}
1-\varepsilon \mu_{1}y_{0}C_{4}\sup_{\mu_{1}X}|x_{s}| \leq
x_{s}(s,t)\leq\varepsilon
\mu_{1}y_{0}C_{4}\sup_{\mu_{1}X}|x_{s}|+1
\end{equation*}
for all $(s,t)\in\mu_{1}X$.  Hence for $\varepsilon$ sufficiently
small, $x_{s}(s,t)>0$ in $\mu_{1}X$.  We have now shown that
$\rho$ is a diffeomorphism.  Moreover, by Lemma 1.1 and the
inverse function theorem $\rho,\rho^{-1}\in C^{r}$.\par
   We now calculate $b_{11}^{4}$ and $b_{1}^{4}$.  We have
\begin{eqnarray}
b_{11}^{4}&=&\frac{(\varepsilon
w_{xy}+\varepsilon^{2n}H^{n}P_{12})^{2}+\varepsilon^{2(n+1)}H^{n+1}P}
{(1+\varepsilon Q)^{2}}\xi_{x}^{2}\\
& &- \frac{2(\varepsilon w_{xy}+\varepsilon^{2n}H^{n}P_{12})}
{1+\varepsilon Q}\xi_{x}\xi_{y} +\xi_{y}^{2}.\nonumber
\end{eqnarray}
Since $\xi_{y}=-b_{12}^{3}\xi_{x}$, plugging into (11) we obtain
\begin{equation*}
b_{11}^{4}=\frac{\varepsilon^{2(n+1)}H^{n+1}P\xi_{x}^{2}}
{(1+\varepsilon Q)^{2}} :=\varepsilon^{2(n+1)}H^{n+1}P_{11}^{4}.
\end{equation*}
To show that $P_{11}^{4}\geq C_{1}$, we now estimate $\xi_{x}$. By
differentiating (10) with respect to $x$, we obtain
\begin{equation*}
b_{12}^{3}(\xi_{x})_{x}+(\xi_{x})_{y}=-(b_{12}^{3})_{x}\xi_{x},
\text{ }\text{ }\text{ }\xi_{x}(x,0)=1.
\end{equation*}
As above let $(x(t),y(t))$ be the parameterization for an
arbitrary characteristic, then $\xi_{x}(t)=\xi_{x}(x(t),y(t))$
satisfies $\dot{\xi_{x}}=-(b_{12}^{3})_{x}\xi_{x}$,
$\xi_{x}(0)=1$.  By the mean value theorem
\begin{equation*}
|\xi_{x}(t)-1|=|\xi_{x}(t)-\xi_{x}(0)|\leq \mu_{1}y_{0}
\sup_{X_{1}}|(b_{12}^{3})_{x}|\sup_{X_{1}} |\xi_{x}|.
\end{equation*}
Therefore
\begin{equation}
1-\varepsilon \mu_{1}y_{0}C_{5}\sup_{X_{1}}|\xi_{x}|
\leq\xi_{x}(t)\leq\varepsilon
\mu_{1}y_{0}C_{5}\sup_{X_{1}}|\xi_{x}|+1.
\end{equation}
Thus for $\varepsilon$ small $\xi_{x}\geq C_{6}>0$, showing that
$P_{11}^{4}\geq C_{1}$ for some $C_{1}>0$ independent of
$\varepsilon$ and $w$.\par
   We now calculate $b_{1}^{4}$.  We have
\begin{equation}
b_{1}^{4}=b_{11}^{3}\xi_{xx}+2b_{12}^{3}\xi_{xy}
+b_{22}^{3}\xi_{yy}+b_{1}^{3}\xi_{x}+b_{2}^{3}\xi_{y}.
\end{equation}
From (10) we obtain
\begin{equation}
\xi_{xy}=-(b_{12}^{3})_{x}\xi_{x}-b_{12}^{3}\xi_{xx},\text{
}\text{ } \xi_{yy}=-(b_{12}^{3})_{y}\xi_{x}-b_{12}^{3}\xi_{xy}.
\end{equation}
Plugging into (13) produces
\begin{eqnarray}
b_{1}^{4}&=&\frac{\varepsilon^{2(n+1)}H^{n+1}P} {(1+\varepsilon
Q)^{2}}\xi_{xx}+b_{1}^{3}\xi_{x}
+b_{2}^{3}\xi_{y}  \nonumber \\
& &+[\partial_{y}(\frac{\varepsilon
w_{xy}+\varepsilon^{2n}H^{n}P_{12}} {1+\varepsilon
Q})-\frac{1}{2}\partial_{x} (\frac{\varepsilon
w_{xy}+\varepsilon^{2n}H^{n}P_{12}}
{1+\varepsilon Q})^{2}]\xi_{x}\\
&=&\varepsilon^{2n}H^{n}Q_{1}+n\varepsilon^{2n}H^{n-1}Q_{2}
\nonumber  \\
& &+[\partial_{y}(\frac{\varepsilon w_{xy}} {1+\varepsilon
Q})-\frac{1}{2}\partial_{x} (\frac{\varepsilon w_{xy}}
{1+\varepsilon Q})^{2}]\xi_{x},\nonumber
\end{eqnarray}
for some $Q_{1}$, $Q_{2}$.  We now calculate the last term of
(15). From (6) we have
\begin{equation}
\frac{-\varepsilon^{2}w_{xy}^{2}} {(1+\varepsilon
Q)^{2}}=\frac{-\varepsilon w_{yy}(1+\varepsilon Q)
+\varepsilon^{2n}H^{n}Q_{3}+\Phi(w)} {(1+\varepsilon Q)^{2}},
\end{equation}
for some $Q_{3}$.  Then plugging (16) into (15), we obtain
\begin{eqnarray*}
& &\partial_{y}(\frac{\varepsilon w_{xy}}{1+\varepsilon
Q})-\frac{1}{2}\partial_{x} (\frac{\varepsilon w_{xy}}
{1+\varepsilon Q})^{2}\\
&=&\partial_{y}(\frac{\varepsilon w_{xy}} {1+\varepsilon Q})
-\frac{1}{2}\partial_{x} (\frac{\varepsilon w_{yy}} {1+\varepsilon
Q})\\
& &+\varepsilon^{2n}H^{n}Q_{4}+n\varepsilon^{2n}H^{n-1}Q_{5}
+\partial_{x}[\frac{\Phi(w)}{2(1+\varepsilon Q)^{2}}]
\end{eqnarray*}
\begin{eqnarray*}
&=&\frac{\varepsilon/2w_{xyy}(1+\varepsilon
w_{xx})-\varepsilon^{2}
w_{xy}w_{xxy}+\varepsilon^{2}/2w_{yy}w_{xxx}} {(1+\varepsilon
Q)^{2}}\\
& &+\varepsilon^{2n}H^{n}Q_{6}+n\varepsilon^{2n}H^{n-1}Q_{7}
+\partial_{x}[\frac{\Phi(w)}{2(1+\varepsilon
Q)^{2}}]\\
&=&\frac{\partial_{x}}{2(1+\varepsilon Q)^{2}}[\varepsilon w_{yy}
(1+\varepsilon w_{xx})-\varepsilon^{2}w_{xy}^{2}]\\
& &+\varepsilon^{2n}H^{n}Q_{6}+n\varepsilon^{2n}H^{n-1}Q_{7}
+\partial_{x}[\frac{\Phi(w)}{2(1+\varepsilon
Q)^{2}}]\\
&=&\frac{\partial_{x}\Phi(w)}{2(1+\varepsilon Q)^{2}}
+\partial_{x}[\frac{\Phi(w)}{2(1+\varepsilon
Q)^{2}}]\\
& &+\varepsilon^{2n}H^{n}Q_{8}+n\varepsilon^{2n}H^{n-1}Q_{9},
\end{eqnarray*}
for some $Q_{4},\ldots,Q_{9}$.  It follows from (15), that
$b_{1}^{4}$ has the desired form.\par
   To complete the proof of Lemma 1.2, we now show that
$\sum|b_{ij}^{4}|_{C^{12}}+|b_{i}^{4}|_{C^{12}}+|b^{4}|_{C^{12}}\leq
C_{2}$, for some constant $C_{2}$ independent of $\varepsilon$ and
$w$.  In view of the fact that $|w|_{C^{16}}\leq C$, this will be
accomplished by showing that $|\xi|_{C^{14}}\leq C_{7}$ for some
$C_{7}$ independent of $\varepsilon$ and $w$.  By (12) we find
that
\begin{equation*}
\sup_{X_{1}}|\xi_{x}|\leq\frac{1} {1-\varepsilon
C_{5}\mu_{1}y_{0}}:=C_{8}.
\end{equation*}
It follows from (10) that
\begin{equation*}
\sup_{X_{1}}|\xi_{y}|\leq C_{9},
\end{equation*}
where $C_{9}$ is independent of $\varepsilon$ and $w$.\par
   We now estimate $\xi_{xx}$.  Differentiate (10) two times with
respect to $x$ to obtain
\begin{equation*}
b_{12}^{3}(\xi_{xx})_{x}+(\xi_{xx})_{y}=-2(b_{12}^{3})_{x}\xi_{xx}
-(b_{12}^{3})_{xx}\xi_{x},\text{ }\text{ }\text{ }\xi_{xx}(x,0)=0.
\end{equation*}
Then the same procedure that yielded (12), produces
\begin{equation*}
\sup_{X_{1}}|\xi_{xx}|\leq \varepsilon
\mu_{1}y_{0}C_{10}\sup_{X_{1}}|\xi_{xx}| +\varepsilon
\mu_{1}y_{0}C_{11}C_{8},
\end{equation*}
implying that
\begin{equation*}
\sup_{X_{1}}|\xi_{xx}|\leq \frac{\varepsilon
\mu_{1}y_{0}C_{11}C_{8}}{1-\varepsilon
\mu_{1}y_{0}C_{10}}:=C_{12}.
\end{equation*}
Furthermore in light of (14), we can use the estimates for
$\xi_{x}$ and $\xi_{xx}$ to estimate $\xi_{xy}$, and then
subsequently $\xi_{yy}$.  Clearly, we can continue this procedure
to yield $|\xi|_{C^{14}}\leq C_{7}$.  q.e.d.\par\medskip
   We now continue defining the sequence of linear operators
$L_{i}(w)$.  To simplify the coefficient of $\partial_{\xi}$ in
$L_{4}(w)$, we remove the portion of $b_{1}^{4}$ involving
$\Phi(w)$ and define
\begin{eqnarray*}
L_{5}(w)&:=&L_{4}(w)-[\partial_{x}(\frac{\Phi(w)}{2(1+\varepsilon
Q)^{2}})+\frac{\partial_{x}\Phi(w)} {2(1+\varepsilon
Q)^{2}}]\xi_{x}\partial_{\xi}\\
&=&\sum_{i,j}b_{ij}^{5}\partial_{x_{i}x_{j}}+\sum_{i}b_{i}^{5}
\partial_{x_{i}}+b^{5},
\end{eqnarray*}
where $x_{1}=\xi$, $x_{2}=\eta$.\par
   To bring $L_{5}(w)$ into the canonical form (4), we shall
need one more change of variables.\medskip

\textbf{Lemma 1.3.}  \textit{For $\varepsilon$ sufficiently small,
there exists a $C^{r}$ diffeomorphism}
\begin{equation*}
\alpha=\alpha(\xi,\eta),\text{   }\beta=H(\xi,\eta),
\end{equation*}
\textit{of a domain $X_{2}\subset\mu_{1}X$ onto $\mu_{2}X$,
$1<\mu_{2}<\mu_{1}$, such that $\mu_{3}X$ properly contains the
image of $\rho^{-1}(X)$} (\textit{where $\rho^{-1}$ is the
diffeomorphism given by Lemma} 1.2),\textit{ for some $\mu_{3}$,
$1<\mu_{3}<\mu_{2}$. In the new variables $(\alpha,\beta)$,
$L_{5}(w)$ is denoted by $L_{6}(w)$ and is given by}
\begin{equation*}
L_{6}(w)=\sum_{i,j}b_{ij}^{6}\partial_{x_{i}x_{j}}+\sum_{i}b_{i}^{6}
\partial_{x_{i}}+b^{6},
\end{equation*}
\textit{where $x_{1}=\alpha$, $x_{2}=\beta$, and}
\begin{eqnarray*}
b_{11}^{6}&=&\varepsilon^{2(n+1)}\beta^{n+1}P_{11}^{6},\\
b_{12}^{6}=b_{21}^{6}&\equiv& 0,\\
b_{22}^{6}&=&P_{22}^{6},\\
b_{1}^{6}&=&\varepsilon^{2n}\beta^{n}P_{1}^{61}+n\varepsilon^{2n}\beta^{n-1}
P_{1}^{62},\\
b_{2}^{6}&=&\varepsilon P_{2}^{61}+n\varepsilon^{2n}\beta^{n-1}
P_{2}^{62},\\
b^{6}&=&\varepsilon^{2n}\beta^{n}P_{3}^{6},
\end{eqnarray*}
\textit{for some $P_{11}^{6}$, $P_{22}^{6}$, $P_{1}^{61}$,
$P_{1}^{62}$, $P_{2}^{61}$, $P_{2}^{62}$, $P_{3}^{6}$, such that
$P_{11}^{6},P_{22}^{6}\geq C_{13}$ for some constant $C_{13}>0$
independent of $\varepsilon$ and $w$. Furthermore
$\sum|b_{ij}^{6}|_{C^{12}}+|b_{i}^{6}|_{C^{12}}+|b^{6}|_{C^{12}}\leq
C_{14}$, for some $C_{14}$ independent of $\varepsilon$ and
$w$.}\medskip

\textit{Proof.}  Using the chain rule we find that
$b_{12}^{6}=b_{11}^{5}\beta_{\xi}\alpha_{\xi}
+b_{22}^{5}\beta_{\eta}\alpha_{\eta}$. Therefore, we seek a smooth
function $\alpha(\xi,\eta)$ such that
\begin{equation}
b_{12}^{6}=b_{11}^{5}\beta_{\xi}\alpha_{\xi}
+b_{22}^{5}\beta_{\eta}\alpha_{\eta}=0\text{ }\text{ }\text{ in
}\text{ }X_{2},\text{ }\text{ }\text{ }\alpha(\xi,0)=\xi,
\end{equation}
where $X_{2}$ will be defined below. By our original assumption on
$H$ made in the introduction, $H_{y}\geq C_{15}$ for some
$C_{15}>0$ independent of $\varepsilon$. Therefore
\begin{equation*}
H_{\eta}=H_{x}\frac{\partial x}{\partial\eta} +H_{y}\frac{
\partial y}{\partial\eta}=-H_{x}\frac{\xi_{y}}{\xi_{x}}+H_{y}
\geq \varepsilon C_{16}+C_{15}\geq C_{17}>0,
\end{equation*}
for some $C_{16}$, $C_{17}$ independent of $\varepsilon$.  Since
$b_{22}^{5}\equiv 1$, it follows that the line $\eta=0$ is
noncharacteristic for (17).  Therefore, the methods used in the
proof of Lemma 1.2 show that the desired function
$\alpha(\xi,\eta)$ exists.\par
   We now define $X_{2}$.  Since
$H_{\eta}\geq C_{17}>0$, we may choose $\mu_{1}>\mu_{2}>1$ such
that the curves $H(\xi,\eta)=\pm\mu_{2}y_{0}$ are properly
contained in the strips $\{(\xi,\eta)\mid
y_{0}\leq\eta\leq\mu_{1}y_{0}\}$, $\{(\xi,\eta)\mid
-y_{0}\geq\eta\geq-\mu_{1}y_{0}\}$.  Then define
$X_{2}\subset\mu_{1}X$ to be the domain in the $\xi,\eta$ plane
bounded by the curves $H(\xi,\eta)=\pm\mu_{2}y_{0}$ and the
characteristic curves of (17) passing through the points
$(\pm\mu_{2}x_{0},0)$. Then the methods of the proof of Lemma 1.2
show that the mapping
$\tau:(\xi,\eta)\mapsto(\alpha(\xi,\eta),\beta(\xi,\eta))$ is a
$C^{r}$ diffeomorphism from $X_{2}$ onto $\mu_{2}X$. Furthermore,
since $\rho^{-1}(X)\subset X_{2}$, if $\mu_{3}$ is chosen large
then $\tau(\rho^{-1}(X))\subset\mu_{3}X$.\par
   We now compute the coefficients $b_{ij}^{6}$, $b_{i}^{6}$,
$b^{6}$.  We have
\begin{eqnarray*}
b_{11}^{6}&=&\varepsilon^{2(n+1)}\beta^{n+1}P_{11}^{4}\alpha_{\xi}^{2}+\alpha_{\eta}^{2}\\
&=&\varepsilon^{2(n+1)}\beta^{n+1}P_{11}^{4}\alpha_{\xi}^{2}
+\varepsilon^{4(n+1)}\beta^{2(n+1)}(P_{11}^{4})^{2}\frac{\beta_{\xi}^{2}}
{\beta_{\eta}^{2}}\alpha_{\xi}^{2}\\
&=&\varepsilon^{2(n+1)}\beta^{n+1}[P_{11}^{4}+
\varepsilon^{2(n+1)}\beta^{n+1}(P_{11}^{4})^{2}\frac{\beta_{\xi}^{2}}
{\beta_{\eta}^{2}}]\alpha_{\xi}^{2}\\
&:=&\varepsilon^{2(n+1)}\beta^{n+1}P_{11}^{6}.
\end{eqnarray*}
As in the proof of Lemma 1.2, $\alpha_{\xi}\geq C_{18}$ for some
$C_{18}>0$ independent of $\varepsilon$ and $w$.  Thus, if
$\varepsilon$ is sufficiently small the properties of $P_{11}^{4}$
imply that $P_{11}^{6}\geq C_{13}$ for some $C_{13}>0$ independent
of $\varepsilon$ and $w$.  Next we calculate $b_{22}^{6}$:
\begin{equation*}
b_{22}^{6}=\varepsilon^{2(n+1)}\beta^{n+1}P_{11}^{4}\beta_{\xi}^{2}
+\beta_{\eta}^{2}:=P_{22}^{6}.
\end{equation*}
Since $H_{\eta}\geq C_{17}$, if $\varepsilon$ is sufficiently
small then $P_{22}^{6}\geq C_{13}$.  Furthermore, by (17)
\begin{eqnarray*}
b_{1}^{6}&=& b_{11}^{5}\alpha_{\xi\xi}+\alpha_{\eta\eta}+b_{1}^{5}
\alpha_{\xi}+b_{2}^{5}\alpha_{\eta}\\
&=& b_{11}^{5}\alpha_{\xi\xi}-\partial_{\eta}(
\frac{\varepsilon^{2(n+1)}\beta^{n+1}P_{11}^{4}\beta_{\xi}\alpha_{\xi}}
{\beta_{\eta}})+b_{1}^{5}\alpha_{\xi}+b_{2}^{5}\alpha_{\eta}\\
&:=&\varepsilon^{2n}\beta^{n}P_{1}^{61}+n\varepsilon^{2n}\beta^{n-1}
P_{1}^{62}.
\end{eqnarray*}
Lastly since $\beta_{\eta}=H_{x}(\frac{-\xi_{y}}{\xi_{x}})+H_{y}=
O(\varepsilon)+H_{y}$, we have
\begin{equation*}
\beta_{\eta\eta}=O(\varepsilon)+H_{yy}=O(\varepsilon)+\varepsilon^{2}\widetilde{H}_{vv}
=O(\varepsilon).
\end{equation*}
Thus
\begin{eqnarray*}
b_{2}^{6}&=&b_{11}^{5}\beta_{\xi\xi}+\beta_{\eta\eta}+b_{1}^{5}
\beta_{\xi}+b_{2}^{5}\beta_{\eta}\\
&:=&\varepsilon P_{2}^{61}+n\varepsilon^{2n}\beta^{n-1}P_{2}^{62}.
\end{eqnarray*}\par
   We complete the proof by noting that the methods of the proof
of Lemma 1.2 show that
$\sum|b_{ij}^{6}|_{C^{12}}+|b_{i}^{6}|_{C^{12}}+|b^{6}|_{C^{12}}\leq
C_{14}$, for some $C_{14}$ independent of $\varepsilon$ and $w$.
q.e.d.\par\medskip
   To obtain the canonical form (4), we define
\begin{eqnarray*}
L_{7}(w)&:=&\frac{1}{b_{22}^{6}}L_{6}(w)\\
&=&\sum_{i,j}b_{ij}^{7}\partial_{x_{i}x_{j}}+\sum_{i}b_{i}^{7}
\partial_{x_{i}}+b^{7},
\end{eqnarray*}
where $x_{1}=\alpha$, $x_{2}=\beta$, and
\begin{eqnarray*}
b_{11}^{7}&=&\varepsilon^{2(n+1)}\beta^{n+1}P_{11}^{7},\\
b_{12}^{7}=b_{21}^{7}&\equiv& 0,\\
b_{22}^{7}&\equiv& 1,\\
b_{1}^{7}&=&\varepsilon^{2n}\beta^{n}P_{1}^{71}+n\varepsilon^{2n}\beta^{n-1}
P_{1}^{72},\\
b_{2}^{7}&=&\varepsilon P_{2}^{71}+n\varepsilon^{2n}\beta^{n-1}
P_{2}^{72},\\
b^{7}&=&\varepsilon^{2n}\beta^{n}P_{3}^{7},
\end{eqnarray*}
for some $P_{11}^{7}$, $P_{1}^{71}$, $P_{1}^{72}$, $P_{2}^{71}$,
$P_{2}^{72}$, $P_{3}^{7}$, such that $P_{11}^{7}\geq C_{19}$ for
some constant $C_{19}>0$ independent of $\varepsilon$ and $w$.  In
the following section, we shall study the existence and regularity
theory for the operator $L_{7}(w)$.
\bigskip
\begin{center}
\textbf{2.  Linear Theory}
\end{center}
\bigskip

   In this section we study the existence and regularity theory for
the operator $L_{7}$.  More precisely, we will first extend the
coefficients of $L_{7}$ onto the entire plane in a manner that
facilitates an a priori estimate, and then prove the existence of
weak solutions having regularity in the $\alpha$-direction.  It
will then be shown that these weak solutions are also regular in
the $\beta$-direction via a boot-strap argument.\par
   For simplicity of notation, put $x=\alpha$, $y=\beta$, and
$\overline{L}=L_{7}(w)$.  Then
\begin{eqnarray*}
\overline{L}&=&\varepsilon^{2(n+1)}y^{n+1}B_{1}\partial_{xx}+\partial_{yy}
+(\varepsilon^{2n}y^{n}B_{2}+n\varepsilon^{2n}y^{n-1}B_{3})\partial_{x}\\
& &+(\varepsilon B_{4}+n\varepsilon^{2n}y^{n-1}B_{5})\partial_{y}
+\varepsilon^{2n}y^{n}B_{6}\\
&:=&\overline{A}\partial_{xx}+\partial_{yy}+\overline{D}\partial_{x}
+\overline{E}\partial_{y}+\overline{F},
\end{eqnarray*}
for some $B_{1},\ldots,B_{6}\in C^{r}$ such that $B_{1}\geq M$ and
$|B_{i}|_{C^{12}}\leq M^{'}$, for some constants $M,M^{'}>0$
independent of $\varepsilon$ and $w$.  By Lemma 1.3
$\overline{A},\overline{D},\overline{E}$, and $\overline{F}$ are
defined in the rectangle $\mu_{2}X$.  We will modify these
coefficients on $\mathbb{R}^{2}-\mu_{2}X$, so that they will be
defined and of class $C^{r}$ on the entire plane.\par
   Choose values $y_{1},\ldots,y_{6}$ such that
$0<y_{1}<\cdots<y_{6}$ and $y_{1}=\mu_{3}y_{0}$,
$y_{6}=\mu_{2}y_{0}$.  Let $\delta, M_{1}>0$ be constants, where
$\delta$ will be chosen small.  Fix a nonnegative cut-off function
$\phi\in C^{\infty}(\mathbb{R})$ such that
$$
\phi(y)=\begin{cases}
1 & \text{if $|y|\leq y_{5}$},\\
0 & \text{if $|y|\geq y_{6}$}.
\end{cases}
$$
Furthermore, define functions $\psi_{1},\psi_{2},\psi_{3}\in
C^{\infty}(\mathbb{R})$ with properties:\bigskip

$i)$ $\psi_{1}(y)=\begin{cases}
0 & \text{if }|y|\leq y_{2},\\
-1 & \text{if }y\leq-y_{3},\\
1 & \text{if }y\geq y_{3},
\end{cases}$\bigskip

$ii)$ $\psi_{1}\leq 0$ if $y\leq 0$, $\psi_{1}\geq 0$ if $y\geq
0$, and $\psi_{1}^{'}\geq 0$,\bigskip

$iii)$ $\psi_{2}(y)=\begin{cases}
0 & \text{if }y\geq-y_{5},\\
-\delta y-\delta(\frac{y_{5}+y_{6}}{2}) & \text{if }y\leq-y_{6},
\end{cases}$\bigskip

$iv)$ $\psi_{2}\geq 0$, and $-\delta\leq\psi_{2}^{'}\leq
0$,\bigskip

$v)$ $\psi_{3}(y)=\begin{cases}
0 & \text{if }|y|\leq y_{3},\\
M_{1} & \text{if }y\leq-y_{4},\\
-M_{1} & \text{if }y\geq y_{4},
\end{cases}$\bigskip

$vi)$ $\psi_{3}\geq 0$ if $y\leq 0$, $\psi_{3}\leq 0$ if $y\geq
0$, and $\psi_{3}^{'}\leq 0$.\bigskip

\noindent Now define smooth extensions of
$\overline{A},\overline{D},\overline{E}$, and $\overline{F}$ to
the entire plane by
\begin{eqnarray*}
A&=&\psi_{1}(y)+\phi(x)\phi(y)\overline{A},\\
D&=&\phi(x)\phi(y)\overline{D},\\
E&=&\psi_{2}(y)+\phi(x)\phi(y)\overline{E},\\
F&=&\psi_{3}(y)+\phi(x)\phi(y)\overline{F},
\end{eqnarray*}
and set
\begin{equation*}
L=A\partial_{xx}+\partial_{yy}+D\partial_{x}+E\partial_{y}+F.
\end{equation*}\par
   Before making estimates for $L$, we must define the function
spaces that will be utilized.  For $m,l\in\mathbb{Z}_{\geq 0}$,
let
\begin{equation*}
C^{(m,\text{
}l)}(\mathbb{R}^{2})=\{u:\mathbb{R}^{2}\rightarrow
\mathbb{R}\mid\partial_{x}^{s}\partial_{y}^{t}u\in
C(\mathbb{R}^{2}),\text{ }s\leq m,\text{ }t\leq l\},
\end{equation*}
and
\begin{equation*}
C^{(m,\text{ }l)}_{c}(\mathbb{R}^{2})=\{u\in C^{(m,\text{
}l)}(\mathbb{R}^{2})\mid u\text{ has compact support}\}.
\end{equation*}
Let $\theta>0$ be a small parameter, and define the norm
\begin{equation*}
\parallel u\parallel_{(m,\text{ }l)}=(\sum_{s\leq m,\text{ }t\leq
l}\theta^{s}\parallel\partial_{x}^{s}\partial_{y}^{t}u
\parallel^{2}_{L^{2}(\mathbb{R}^{2})})^{1/2}.
\end{equation*}
Then define $H_{\theta}^{(m,\text{ }l)}(\mathbb{R}^{2})$ to be the
closure of $C_{c}^{(m,\text{ }l)}(\mathbb{R}^{2})$ in the norm
$\parallel\cdot\parallel_{(m,\text{ }l)}$.  Furthermore, let
$H^{m}(\mathbb{R}^{2})$ be the Sobolev space with square
integrable derivatives up to and including order $m$, with norm
$\parallel\cdot\parallel_{m}$.  Lastly, denote the
$L^{2}(\mathbb{R}^{2})$ inner product and norm by $(\cdot,\cdot)$
and $\parallel\cdot\parallel$ respectively.\par
   We are now ready to establish a basic estimate for the
operator $L$ on $\mathbb{R}^{2}$.  This estimate will be used to
establish a more general estimate, which will in turn be used as
the foundation for the proof of the existence of weak
solutions.\medskip

\textbf{Lemma 2.1.}  \textit{If $\varepsilon$ is sufficiently
small, then there exists a constant $C_{1}>0$ independent of
$\varepsilon$, and functions $a(y),b(y),\gamma(y)\in
C^{\infty}(\mathbb{R})$ where $\gamma=O(1)$ as
$y\rightarrow\infty$, and $\gamma=O(|y|)$ as $y\rightarrow-\infty$
such that}
\begin{equation*}
(au+bu_{y},Lu)\geq C_{1}(\parallel\gamma
u_{y}\parallel^{2}+\parallel u\parallel^{2}),\textit{ for all
$u\in C_{c}^{\infty}(\mathbb{R}^{2})$.}
\end{equation*}

\textit{Proof.}  We first define the functions $a$ and $b$.  Let
$M_{2},M_{3},M_{4}>0$ be constants satisfying $M_{3}<M_{2}$ and
$\frac{1}{2}M_{4}-M_{2}\geq 1$.  Then choose $a,b\in
C^{\infty}(\mathbb{R})$ and $M_{2},M_{3},M_{4}$ such that:\bigskip

$i)$ $a(y)=\begin{cases}
y^{2}-M_{2} & \text{if }|y|\leq y_{5},\\
-M_{3} & \text{if }|y|\geq y_{6},
\end{cases}$\bigskip

$ii)$ $a\leq-M_{3}$, $a^{'}\leq 0$ if $y\leq 0$, $a^{'}\geq 0$ if
$y\geq 0$, and $a^{''}\geq-\delta$,\bigskip

$iii)$ $b(y)=\begin{cases}
1 & \text{if }y\geq 0,\\
-M_{4}y+1 & \text{if }y\leq-y_{2},
\end{cases}$\bigskip

$iv)$ $b\geq 1$, and $b^{'}\leq 0$.\bigskip

\noindent Now let $u\in C_{c}^{\infty}(\mathbb{R}^{2})$, and
integrate by parts to obtain
\begin{equation*}
(au+bu_{y},Lu)=\int\int_{\mathbb{R}^{2}}I_{1}u_{x}^{2}+2I_{2}u_{x}u_{y}
+I_{3}u_{y}^{2}+I_{4}u^{2},
\end{equation*}
where
\begin{eqnarray*}
I_{1}&=&(\frac{1}{2}b^{'}-a)A+\frac{1}{2}bA_{y},\\
I_{2}&=&-\frac{1}{2}bA_{x}+\frac{1}{2}bD,\\
I_{3}&=&-a-\frac{1}{2}b^{'}+bE,\\
I_{4}&=&\frac{1}{2}aA_{xx}
+\frac{1}{2}a^{''}-\frac{1}{2}aD_{x}-\frac{1}{2}(aE)_{y}-(\frac{1}{2}b^{'}-a)F
-\frac{1}{2}bF_{y}.
\end{eqnarray*}\par
   We now estimate $I_{1}$.  If $|y|\leq y_{3}$ then
\begin{eqnarray*}
I_{1}&\geq&[(M_{2}-y^{2})\varepsilon^{2(n+1)}y^{n+1}B_{1}
+\frac{(n+1)}{2}\varepsilon^{2(n+1)}y^{n}B_{1}\\
& &+\frac{1}{2}\varepsilon^{2(n+1)}y^{n+1}b\partial_{y}B_{1}]\phi(x)\\
&=&\varepsilon^{2(n+1)}y^{n}[(M_{2}-y^{2})yB_{1}+\frac{(n+1)}{2}B_{1}+
\frac{1}{2}yb\partial_{y}B_{1}]\phi(x)\\
&\geq&\varepsilon^{2(n+1)}C_{2}y^{n}\phi(x)\geq 0,
\end{eqnarray*}
for some constants $C_{2}>0$, if $y_{3}$ is chosen sufficiently
small.  Moreover, if $|y|\geq y_{3}$ we have
\begin{equation*}
I_{1}\geq O(\varepsilon^{2(n+1)})+\begin{cases}
M_{3} & \text{if }y\geq 0\\
\frac{1}{2}M_{4}-M_{2} & \text{if }y<0
\end{cases}
\geq C_{3},
\end{equation*}
for some $C_{3}>0$, if $\varepsilon$ is small.\par
   To estimate $I_{3}$, we observe that for $|y|\leq y_{6}$,
\begin{equation*}
I_{3}\geq M_{3}+O(\varepsilon).
\end{equation*}
Furthermore, if $|y|\geq y_{6}$ then
\begin{equation*}
I_{3}\geq M_{3}+\begin{cases}
0 & \text{if }y\geq 0,\\
\delta M_{4}y^{2} & \text{if }y<0.
\end{cases}
\end{equation*}
Hence, $I_{3}\geq\gamma^{2}(y)$ for some $\gamma\in
C^{\infty}(\mathbb{R})$ such that $\gamma=O(1)$ as
$y\rightarrow\infty$, and $\gamma=O(|y|)$ as
$y\rightarrow-\infty$.\par
   Next we show that
\begin{equation*}
\int\int_{\mathbb{R}^{2}}I_{1}u_{x}^{2}+2I_{2}u_{x}u_{y}
+I_{3}u_{y}^{2}\geq C_{4}\parallel\gamma u_{y}\parallel^{2},
\end{equation*}
for some $C_{4}>0$.  From our estimates on $I_{1}$ and $I_{3}$,
this will follow if $I_{1}I_{3}-2I_{2}^{2}\geq 0$.  A calculation
shows that when $|y|\leq y_{6}$, we have
\begin{eqnarray*}
I_{1}I_{3}-2I_{2}^{2}&\geq&\varepsilon^{2(n+1)}C_{5}y^{n}\phi(x)
+O(n\varepsilon^{2n}y^{n-1}\phi(x)+\varepsilon^{2n}y^{n}|\phi^{'}(x)|)^{2}\\
&=&\varepsilon^{2(n+1)}y^{n}[C_{5}+\varepsilon^{2n-2}O(n^{2}y^{n-2}\phi(x)
+y^{n}|\phi^{'}(x)|^{2}\phi^{-1}(x)\\
& &+ny^{n-1}|\phi^{'}(x)|)]\phi(x)\\
&\geq& 0,
\end{eqnarray*}
for some $C_{5}>0$ independent of $\varepsilon$, if $\varepsilon$
is sufficiently small.  Moreover, if $|y|\geq y_{6}$ then
\begin{equation*}
I_{1}I_{3}-2I_{2}^{2}=I_{1}I_{3}>0,
\end{equation*}
from which we obtain the desired conclusion.\par
   Lastly, we estimate $I_{4}$.  In the strip $|y|\leq y_{4}$, we
obtain
\begin{equation*}
I_{4}\geq 1+O(\varepsilon).
\end{equation*}
Furthermore, if $|y|\geq y_{4}$ then
\begin{equation*}
I_{4}\geq\begin{cases}
M_{1}M_{3}+O(\varepsilon+\delta) & \text{if }y\geq 0,\\
M_{1}(\frac{1}{2}M_{4}-M_{2})+O(\varepsilon+\delta) & \text{if
}y<0.
\end{cases}
\end{equation*}
Therefore, $I_{4}\geq C_{6}$ for some $C_{6}>0$ independent of
$\varepsilon$.  q.e.d.\par\medskip
   Having established the basic estimate, our goal shall now be to
establish a more general estimate that involves derivatives of
higher order in the $x$-direction.  Let $<\cdot,\cdot>_{m}$ denote
the inner product on $H_{\theta}^{(m,0)}(\mathbb{R}^{2})$, that
is,
\begin{equation*}
<u,v>_{m}=\int\int_{\mathbb{R}^{2}}\sum_{s=0}^{m}\theta^{s}
\partial_{x}^{s}u\partial_{x}^{s}v,\text{ for all }u,v\in
H_{\theta}^{(m,0)}(\mathbb{R}^{2}).
\end{equation*}

\textbf{Theorem 2.1.}  \textit{If $\varepsilon=\varepsilon(m)$ is
sufficiently small, then for each $m\leq r-2$, there exist
constants $\theta(m)>0$ and $C_{m}>0$, both depending on
$|A|_{C^{m+2}(\mathbb{R}^{2})},|D|_{C^{m+2}(\mathbb{R}^{2})}$,
$|E|_{C^{m+2}(\mathbb{R}^{2})}$, and
$|F|_{C^{m+2}(\mathbb{R}^{2})}$, such that for all
$\theta\leq\theta(m)$}
\begin{equation*}
<au+bu_{y},Lu>_{m}\geq C_{m}(\parallel u\parallel_{(m,0)}^{2}
+\sum_{s=0}^{m}\theta^{s}\parallel\gamma\partial_{x}^{s}u_{y}\parallel^{2}),
\end{equation*}
\textit{for all $u\in C_{c}^{\infty}(\mathbb{R}^{2})$}.\medskip

\textit{Proof.}  We shall prove the estimate by induction on $m$.
The case $m=0$ is given by Lemma 2.1.  Let $m\geq 1$, and assume
that the estimate holds for all integers less than $m$.\par
   Let $u\in C_{c}^{\infty}(\mathbb{R}^{2})$ and set
$w=\partial_{x}^{m}u$, then
\begin{eqnarray}
& &<au+bu_{y},Lu>_{m}\\
&=&<au+bu_{y},Lu>_{m-1}+\theta^{m}(aw+bw_{y},L_{m}w)\nonumber\\
& &+\theta^{m}(a\partial_{x}^{m}u+b\partial_{x}^{m}u_{y},
\sum_{i=0}^{m-1}\partial_{x}^{i}(E_{x}\partial_{x}^{m-1-i}u_{y}+
\partial_{x}F_{m-1-i}\partial_{x}^{m-1-i}u)),\nonumber
\end{eqnarray}
where
\begin{eqnarray*}
L_{m}&=&A\partial_{xx}+\partial_{yy}+D_{m}\partial_{x}
+E\partial_{y}+F_{m},\\
D_{m}&=&D+mA_{x},\text{ }\text{ }\text{ }F_{m}=F+
mD_{x}+\frac{m(m-1)}{2}A_{xx}.
\end{eqnarray*}
We now estimate each term on the right-hand side of (18).  By the
induction assumption,
\begin{equation}
<au+bu_{y},Lu>_{m-1}\geq C_{m-1}(\parallel
u\parallel_{(m-1,0)}^{2}+\sum_{s=0}^{m-1}\theta^{s}\parallel
\gamma\partial_{x}^{s}u_{y}\parallel^{2}).
\end{equation}
In addition, since $D_{x},A_{x},A_{xx}$ have compact support and
both $D_{m}=O(mn\varepsilon^{2n}y^{n-1})$, and
$mD_{x}+\frac{m(m-1)}{2}A_{xx}=O(m^{2}n\varepsilon^{2n})$ near the
origin, if $\varepsilon=\varepsilon(m)$ is sufficiently small then
the coefficients of $L_{m}$ have the same properties as those of
$L$ so that Lemma 2.1 applies to yield,
\begin{equation}
\theta^{m}(aw+bw_{y},L_{m}w)\geq\theta^{m} C_{1}(\parallel\gamma
w_{y}\parallel^{2}+\parallel w\parallel^{2}).
\end{equation}
Furthermore, integrating by parts produces
\begin{eqnarray}
& &(a\partial_{x}^{m}u+b\partial_{x}^{m}u_{y},
\sum_{i=0}^{m-1}\partial_{x}^{i}(E_{x}\partial_{x}^{m-1-i}u_{y}+
\partial_{x}F_{m-1-i}\partial_{x}^{m-1-i}u))\\
&=& \int\int_{\mathbb{R}^{2}}[e_{m-1}(\partial_{x}^{m-1}u)^{2}+
e_{m-2}(\partial_{x}^{m-2}u)^{2}+\cdots+e_{0}u^{2}\nonumber\\
& &+f_{m-1}(\partial_{x}^{m-1}u_{y})^{2}+
f_{m-2}(\partial_{x}^{m-2}u_{y})^{2}+\cdots+f_{0}u_{y}^{2}\nonumber\\
& &+g_{m-1}\partial_{x}^{m}u\partial_{x}^{m-1}u_{y}
+g_{m-2}\partial_{x}^{m-1}u\partial_{x}^{m-2}u_{y}+\cdots+
g_{0}u_{x}u_{y}],\nonumber
\end{eqnarray}
for some functions $e_{i},f_{i},g_{i}$ depending on the
derivatives of $A,D,E$ and $F$ up to and including order
$m+2$.\par
   Observe that the power of $\theta$ in the third term on the
right of (18), is sufficiently large to guarantee that the
right-hand side of (21) may be absorbed into the combined
right-hand sides of (19) and (20), for all $\theta<\theta(m)$ if
$\theta(m)$ is sufficiently small.  Thus, we obtain
\begin{equation*}
<au+bu_{y},Lu>_{m}\geq C_{m}(\parallel u\parallel_{(m,0)}^{2}
+\sum_{s=0}^{m}\theta^{s}\parallel\gamma\partial_{x}^{s}u_{y}\parallel^{2}),
\end{equation*}
completing the proof by induction.  q.e.d.\par\medskip
   Let $f\in L^{2}(\mathbb{R}^{2})$, and consider the equation
\begin{equation}
Lu=f.
\end{equation}
A function $u\in L^{2}(\mathbb{R}^{2})$ is said to be a weak
solution of (22) if
\begin{equation*}
(u,L^{*}v)=(f,v),\text{ for all }v\in
C_{c}^{\infty}(\mathbb{R}^{2}),
\end{equation*}
where $L^{*}$ is the formal adjoint of $L$.  The estimate of
Theorem 2.1 shall serve as the basis for establishing the
existence of weak solutions via the method of Galerkin
approximation.  That is, we shall construct certain
finite-dimensional approximations of (22), and then pass to the
limit to obtain a solution.\par
   Let $\{\phi_{l}\}_{l=1}^{\infty}$ be a basis of
$H_{\theta}^{2m+2}(\mathbb{R})$ that is orthonormal in
$H_{\theta}^{m}(\mathbb{R})$.  Such a sequence may be constructed
by applying the Gram-Schmidt process to a basis of
$H_{\theta}^{2m+2}(\mathbb{R})$.  Choose a positive integer $N$.
We seek an approximate solution, $u^{N}$, of equation (22) in the
form
\begin{equation*}
u^{N}(x,y)=\sum_{l=1}^{N}d_{l}^{N}(y)\phi_{l}(x),
\end{equation*}
where the functions $d_{l}^{N}$ are to be determined from the
relations
\begin{equation}
\int_{\mathbb{R}}\sum_{s=0}^{m}\theta^{s}\frac{d^{s}\phi_{l}}{dx^{s}}\partial_{x}^{s}
Lu^{N}dx=\int_{\mathbb{R}}\sum_{s=0}^{m}\theta^{s}\frac{d^{s}\phi_{l}}{dx^{s}}
\partial_{x}^{s}fdx,\text{ }\text{ }\text{ }l=1,\ldots,N.
\end{equation}
The following lemma will establish the existence of the
$d_{l}^{N}$.\medskip

\textbf{Lemma 2.2.}  \textit{Suppose that
$\varepsilon=\varepsilon(m)$ and $\theta(m)$ are sufficiently
small, and $f\in H_{\theta}^{(m,0)}(\mathbb{R}^{2})$, $m\leq r-2$.
Then there exist functions $d_{l}^{N}\in H^{2}(\mathbb{R})$,
$l=1,\ldots,N$, satisfying} (23) \textit{in the
$L^{2}(\mathbb{R})$-sense}.\medskip

\textit{Proof.}  Choose $\varepsilon$ and $\theta$ so small that
Theorem 2.1 is valid.  Since $\{\phi_{l}\}_{l=1}^{\infty}$ is an
orthonormal set in $H_{\theta}^{m}(\mathbb{R})$, (23) becomes
\begin{eqnarray}
& &(d_{l}^{N})^{''}+\sum_{i=1}^{N}\sum_{s=0}^{m}
(\int_{\mathbb{R}}\theta^{s}
\frac{d^{s}\phi_{l}}{dx^{s}}\partial_{x}^{s}(E\phi_{i})dx)(d_{i}^{N})^{'}\nonumber\\
& &+\sum_{i=1}^{N}\sum_{s=0}^{m} \theta^{s}(\int_{\mathbb{R}}
\frac{d^{s}\phi_{l}}{dx^{s}}\partial_{x}^{s}(A\phi_{i}^{''})
+\frac{d^{s}\phi_{l}}{dx^{s}}\partial_{x}^{s}(D\phi_{i}^{'})
+\frac{d^{s}\phi_{l}}{dx^{s}}\partial_{x}^{s}(F\phi_{i})dx)d_{i}^{N}\\
&
&=\int_{\mathbb{R}}\sum_{s=0}^{m}\theta^{s}\frac{d^{s}\phi_{l}}{dx^{s}}
\partial_{x}^{s}fdx,\text{ }\text{ }\text{ }l=1,\ldots,N.\nonumber
\end{eqnarray}
By the theory of ordinary differential equations, it is sufficient
to prove uniqueness to obtain existence of a solution to system
(24).\par
   We now establish the uniqueness of solutions to (24) in the
space $H^{2}(\mathbb{R})$.  Multiply (23) by
$a(y)d_{l}^{N}(y)+b(y)(d_{l}^{N})^{'}(y)$, sum over $l$ from $1$
to $N$, and then integrate with respect to $y$ over $\mathbb{R}$
to obtain
\begin{equation*}
<au^{N}+bu_{y}^{N},Lu^{N}>_{m}=<au^{N}+bu_{y}^{N},f>_{m}.
\end{equation*}
It now follows from Theorem 2.1 that
\begin{equation}
C_{m}(\parallel u^{N}\parallel_{(m,0)}^{2}+
\sum_{s=0}^{m}\theta^{s}\parallel\gamma\partial_{x}^{s}u_{y}^{N}\parallel^{2})
\leq<au^{N}+bu_{y}^{N},f>_{m},
\end{equation}
for some constant $C_{m}>0$ independent of $N$.  Again using the
orthonormal properties of $\{\phi_{l}\}_{l=1}^{\infty}$, we find
\begin{equation}
\sum_{l=1}^{N}(\parallel d_{l}^{N}\parallel_{\mathbb{R}}^{2}
+\parallel\gamma(d_{l}^{N})^{'}\parallel_{\mathbb{R}}^{2})=\parallel
u^{N}\parallel_{(m,0)}^{2}+
\sum_{s=0}^{m}\theta^{s}\parallel\gamma\partial_{x}^{s}u_{y}^{N}\parallel^{2}.
\end{equation}
Uniqueness for solutions of (23) in the space of functions for
which the left-hand side of (26) is finite, now follows from (25)
and (26).  Thus, existence of a solution in this space is
guaranteed; furthermore, since we can solve for $(d_{l}^{N})^{''}$
in (24), it follows that this solution is in $H^{2}(\mathbb{R})$.
q.e.d.\par\medskip
   Before proving the existence of a weak solution to equation
(22), we will need one more lemma.\medskip

\textbf{Lemma 2.3.}  \textit{Let $v\in
C_{c}^{\infty}(\mathbb{R}^{2})$.  Then there exists a unique
solution, $\widehat{v}\in H^{(\infty,0)}(\mathbb{R}^{2})\cap
C^{\infty}(\mathbb{R}^{2})$, of the equation}
\begin{equation}
(-\theta)^{m}\partial_{x}^{2m}\widehat{v}+(-\theta)^{m-1}
\partial_{x}^{2(m-1)}\widehat{v}+\cdots+\widehat{v}=v.
\end{equation}

\textit{Proof.}  By the Riesz Representation Theorem, there exists
a unique $\widehat{v}\in H^{(m,0)}(\mathbb{R}^{2})$, such that
\begin{equation}
<\widehat{v},w>_{m}=(v,w),\text{ for all }w\in
C_{c}^{\infty}(\mathbb{R}^{2}).
\end{equation}
Thus $\widehat{v}$ is a weak solution of (27), and by the theory
of ordinary differential equations with parameter, we have
$\widehat{v}\in C^{\infty}(\mathbb{R}^{2})$.\par
   We now show that $\widehat{v}\in
H^{(\infty,0)}(\mathbb{R}^{2})$.  It follows from (28) and the
result of Friedrichs [2] on the identity of weak and strong
solutions, that there exists a sequence
$\{\widehat{v}^{k}\}_{k=1}^{\infty}\subset
C_{c}^{\infty}(\mathbb{R}^{2})$ such that
$\widehat{v}^{k}\rightarrow\widehat{v}$ in
$H^{(m,0)}(\mathbb{R}^{2})$, and
\begin{equation*}
(-\theta)^{m}\partial_{x}^{2m}\widehat{v}^{k}+\cdots+
(-\theta)^{\frac{m_{0}+2}{2}}\partial_{x}^{m_{0}+2}\widehat{v}^{k}\rightarrow
v-(-\theta)^{\frac{m_{0}}{2}}\partial_{x}^{m_{0}}\widehat{v}-\cdots-
\widehat{v}
\end{equation*}
in $L^{2}(\mathbb{R}^{2})$, where $m_{0}=m$ if $m$ is even, and
$m_{0}=m-1$ if $m$ is odd.  Therefore
\begin{eqnarray*}
&
&\int\!\!\int_{\mathbb{R}^{2}}v^{2}\\
&=&\!\!\!\!\int\!\!\int_{\mathbb{R}^{2}}
[(-\theta)^{m}\partial_{x}^{2m}\widehat{v}+\cdots+
(-\theta)^{\frac{m_{0}+2}{2}}\partial_{x}^{m_{0}+2}\widehat{v}]^{2}\\
\!\!\!\!&+
&\!\!\!\!\int\!\!\int_{\mathbb{R}^{2}}2[(-\theta)^{m}\partial_{x}^{2m}\widehat{v}+\cdots+
(-\theta)^{\frac{m_{0}+2}{2}}\partial_{x}^{m_{0}+2}\widehat{v}]\\
& &\text{ }\text{ }\text{ }\text{ }\text{
}\cdot[(-\theta)^{\frac{m_{0}}{2}}\partial_{x}^{m_{0}}\widehat{v}+\cdots+
\widehat{v}]\\
\!\!\!\!&+
&\!\!\!\!\int\!\!\int_{\mathbb{R}^{2}}[(-\theta)^{\frac{m_{0}}{2}}\partial_{x}^{m_{0}}\widehat{v}+\cdots+
\widehat{v}]^{2}\\
\!\!\!\!&\geq&\!\!\!\!\lim_{k\rightarrow\infty}\int\!\!\int_{\mathbb{R}^{2}}2
[(-\theta)^{m}\partial_{x}^{2m}\widehat{v}^{k}+\cdots+
(-\theta)^{\frac{m_{0}+2}{2}}\partial_{x}^{m_{0}+2}\widehat{v}^{k}]\\
& &\text{ }\text{ }\text{ }\text{ }\text{ }\text{ }\text{ }\text{
}\text{ }\text{ }\text{
}\cdot[(-\theta)^{\frac{m_{0}}{2}}\partial_{x}^{m_{0}}\widehat{v}^{k}+\cdots+
\widehat{v}^{k}]\\
\!\!\!\!&+ &\!\!\!\!\lim_{k\rightarrow\infty}
\int\!\!\int_{\mathbb{R}^{2}}[(-\theta)^{\frac{m_{0}}{2}}\partial_{x}^{m_{0}}\widehat{v}^{k}+\cdots+
\widehat{v}^{k}]^{2}.
\end{eqnarray*}
Integrating by parts yields
\begin{equation*}
\int\!\!\int_{\mathbb{R}^{2}}\!\!v^{2}\geq
\lim_{k\rightarrow\infty}\int\!\!\int_{\mathbb{R}^{2}}
\theta^{m+1}(\partial_{x}^{m+1}\widehat{v}^{k})^{2}+ \cdots+
(\widehat{v}^{k})^{2},
\end{equation*}
if $m>1$.  Since bounded sets in Hilbert spaces are weakly
compact, $\widehat{v}^{k_{l}}\rightharpoonup\overline{v}$ weakly
in $H^{(m+1,0)}(\mathbb{R}^{2})$, for some $\overline{v}\in
H^{(m+1,0)}(\mathbb{R}^{2})$, where
$\{\widehat{v}^{k_{l}}\}_{l=1}^{\infty}$ is a subsequence of
$\{\widehat{v}^{k}\}$.  For simplicity, we denote
$\widehat{v}^{k_{l}}$ by $\widehat{v}^{k}$.\par
   We now show that $\widehat{v}\equiv\overline{v}$.  By the Riesz
Representation Theorem, there exists $w\in
H^{(m+1,0)}(\mathbb{R}^{2})$ such that
\begin{equation*}
<w,z>_{m+1}=<\widehat{v}-\overline{v},z>_{m},\text{ for all } z\in
H^{(m+1,0)}(\mathbb{R}^{2}).
\end{equation*}
In particular, setting $z=\widehat{v}^{k}-\overline{v}$ we have
\begin{equation}
\lim_{k\rightarrow\infty}<w,\widehat{v}^{k}-\overline{v}>_{m+1}
=\lim_{k\rightarrow\infty}<\widehat{v}-\overline{v},\widehat{v}^{k}
-\overline{v}>_{m}=\parallel\widehat{v}-\overline{v}\parallel_{(m,0)}^{2}.
\end{equation}
Furthermore, since $\widehat{v}^{k}\rightharpoonup\overline{v}$ we
have
\begin{equation}
\lim_{k\rightarrow\infty}<w,\widehat{v}^{k}-\overline{v}>_{m+1}=0.
\end{equation}
Combining (29) and (30) we obtain $\widehat{v}\equiv\overline{v}$
in $H^{(m,0)}(\mathbb{R}^{2})$, implying that $\widehat{v}\in
H^{(m+1,0)}(\mathbb{R}^{2})$. Recall that we assumed that $m>1$;
however, if $m=1$ we still obtain $\widehat{v}\in
H^{(m+1,0)}(\mathbb{R}^{2})$ by solving for
$\partial_{xx}\widehat{v}$ in (27).  A boot-strap argument can now
be used to show that $\widehat{v}\in
H^{(\infty,0)}(\mathbb{R}^{2})$.  q.e.d.\par\medskip
   We are now ready to establish the existence of a weak solution
of equation (22), having regularity in the $x$-direction.\medskip

\textbf{Theorem 2.2.}  \textit{If $\varepsilon=\varepsilon(m)$ and
$\theta(m)$ are sufficiently small, then for every $f\in
H^{(m,0)}_{\theta}(\mathbb{R}^{2})$, $m\leq r-2$, there exists a
unique weak solution $u\in H^{(m,1)}_{\theta}(\mathbb{R}^{2})$ of}
(22).
\medskip

\textit{Proof.}  For each $N\in\mathbb{Z}_{>0}$, let $u^{N}\in
H_{\theta}^{(m,2)}(\mathbb{R}^{2})$ be given by Lemma 2.2.  Then
applying Cauchy's inequality ($pq\leq\kappa
p^{2}+\frac{1}{4\kappa}q^{2}$, $\kappa>0$) to the right-hand side
of (25), we obtain
\begin{equation}
\parallel u^{N}\parallel_{(m,1)}\leq C_{m}^{'}\parallel
f\parallel_{(m,0)},
\end{equation}
where $C_{m}^{'}$ is independent of $N$.  Since bounded sets in
Hilbert spaces are weakly compact, there exists a subsequence
$\{u^{N_{i}}\}_{i=1}^{\infty}$ such that $u^{N_{i}}\rightharpoonup
u$ in $H_{\theta}^{(m,1)}(\mathbb{R}^{2})$, for some $u\in
H_{\theta}^{(m,1)}(\mathbb{R}^{2})$.\par
   We now show that $u$ is a weak solution of (22).  Let $v\in
C_{c}^{\infty}(\mathbb{R}^{2})$ and let $\widehat{v}\in
H^{(\infty,0)}(\mathbb{R}^{2})\cap C^{\infty}(\mathbb{R}^{2})$ be
the solution of
\begin{equation*}
(-\theta)^{m}\partial_{x}^{2m}\widehat{v}+(-\theta)^{m-1}
\partial_{x}^{2(m-1)}\widehat{v}+\cdots+\widehat{v}=v,
\end{equation*}
given by Lemma 2.3.  Since $\{\phi_{l}(x)\}_{l=1}^{\infty}$ forms
a basis in $H_{\theta}^{2m+2}(\mathbb{R})$, we can find
$e_{l}^{N_{*}}(y)\in H^{\infty}(\mathbb{R})$ such that
$v^{N_{*}}:=\sum_{l=1}^{N_{*}}e_{l}^{N_{*}}(y)\phi_{l}(x)
\rightarrow\widehat{v}$ in $H_{\theta}^{(2m+2,2)}(\mathbb{R}^{2})$
as $N_{*}\rightarrow\infty$.  Then multiply (23) by
$e_{l}^{N_{*}}$, sum over $l$ from 1 to $N_{*}$, and integrate
with respect to $y$ over $\mathbb{R}$ to obtain
\begin{equation*}
<v^{N_{*}},Lu^{N_{i}}>_{m}=<v^{N_{*}},f>_{m}.
\end{equation*}
Integrating by parts and letting $N_{i}\rightarrow\infty$
produces,
\begin{equation*}
(u,L^{*}(v^{N_{*}}+\cdots+(-\theta)^{m}\partial_{x}^{2m}v^{N_{*}}))=
(f,v^{N_{*}}+\cdots+(-\theta)^{m}\partial_{x}^{2m}v^{N_{*}}).
\end{equation*}
Furthermore, by letting $N_{*}\rightarrow\infty$ we obtain
\begin{equation*}
(u,L^{*}v)=(f,v).
\end{equation*}
Uniqueness of the weak solution follows from (31).
q.e.d.\par\medskip
   We now prove regularity in the $y$-direction for the weak
solution given by Theorem 2.2, in the case that $f\in
H^{m}(\mathbb{R}^{2})$.  The following standard lemma concerning
difference quotients will be needed.\medskip

\textbf{Lemma 2.4.}  \textit{Let $w\in L^{2}(\mathbb{R}^{2})$ have
compact support, and define}
\begin{equation*}
w^{h}=\frac{1}{h}(w(x,y+h)-w(x,y)).
\end{equation*}
\textit{If $\parallel w^{h}\parallel\leq C_{8}$ where $C_{8}$ is
independent of $h$, then $w\in H^{(0,1)}(V)$ for any compact
$V\subset\mathbb{R}^{2}$. Furthermore, if $w\in
H^{(0,1)}(\mathbb{R}^{2})$ then}
\begin{equation*}
 \parallel w^{h}\parallel\leq
C_{9}\parallel w_{y}\parallel,
\end{equation*}
\textit{for some $C_{9}$ independent of $h$.}\medskip

\textbf{Theorem 2.3.}  \textit{Suppose that the hypotheses of
Theorem} 2.2 \textit{are fulfilled and that $f\in
H^{m}(\mathbb{R}^{2})$, then $u\in H^{m}(\mu_{2}X)$.}\medskip

\textit{Proof.}  If $m=0,1$, then the desired conclusion follows
directly from Theorem 2.2, so assume that $m\geq 2$.  Let
$\zeta\in C^{\infty}(\mathbb{R}^{2})$ be a cut-off function such
that
$$
\zeta(x,y)=\begin{cases}
1 & \text{if $(x,y)\in\mu_{2}X$},\\
0 & \text{if $(x,y)\in(\mu_{2}+1)X$}.
\end{cases}
$$
Let $u\in H_{\theta}^{(m,1)}(\mathbb{R}^{2})$ be the weak solution
of (22) given by Theorem 2.2.  Set $w=\zeta u$, then since $u$ is
a weak solution of (22) we obtain
\begin{equation*}
[w,v]:=\int\int_{\mathbb{R}^{2}}w_{y}v_{y}-Ew_{y}v-Fwv =\int\int
_{\mathbb{R}^{2}}-\widetilde{f}v,\text{ for all }v\in
C_{c}^{\infty}(\mathbb{R}^{2}),
\end{equation*}
where $\widetilde{f}=\zeta f-A\zeta u_{xx}
+\zeta_{yy}u+2\zeta_{y}u_{y}-D\zeta u_{x}+E\zeta_{y}u$.\par
   Using Lemma 2.4 and the fact that $\widetilde{f}\in
L^{2}(\mathbb{R}^{2})$, we have
\begin{eqnarray}
|[w^{h},v]|&\leq&|[w,v^{-h}]|+C_{10}\parallel
v\parallel_{(0,1)}\nonumber\\
&=&|\int\int_{\mathbb{R}^{2}}\widetilde{f}v^{-h}|
+C_{10}\parallel v\parallel_{(0,1)}\\
&\leq&C_{11}\parallel v\parallel_{(0,1)},\nonumber
\end{eqnarray}
for some constants $C_{10}$, $C_{11}$ independent of $h$.
Furthermore, integrating by parts yields
\begin{equation}
C_{12}\parallel v\parallel_{(0,1)}^{2}\leq|[v,v]|+C_{13}\parallel
v\parallel.
\end{equation}
The estimates (32) and (33) also hold if $v=w^{h}$. Therefore
\begin{eqnarray*}
C_{12}\parallel w^{h}\parallel_{(0,1)}^{2}&\leq&C_{11}
\parallel w^{h}\parallel_{(0,1)}+C_{13}\parallel
w^{h}\parallel\\
&\leq&C_{11}\parallel w^{h}\parallel_{(0,1)}+C_{14},
\end{eqnarray*}
for some constant $C_{14}$ independent of $h$.  It follows that
$\parallel w^{h}\parallel_{(0,1)}\leq C_{15}$ independent of $h$.
Hence, by Lemma 2.4 $w\in H^{(0,2)}(V)$ for any compact
$V\subset\mathbb{R}^{2}$.  Since $w\equiv u$ in $\mu_{2}X$, we
have $u\in H^{(0,2)}(\mu_{2}X)$.\par
   By differentiating $Lu=f$ with respect to $x$, $s=1,\ldots,m-2$
times, we obtain
\begin{equation}
L_{s}z=\partial_{x}^{s}f-\sum_{i=0}^{s-1}\partial_{x}^{i}(
E_{x}\partial_{x}^{s-1-i}u_{y}+\partial_{x}F_{s-1-i}\partial_{x}
^{s-1-i}u),
\end{equation}
where $z=\partial_{x}^{s}u$ and $L_{s}$, $F_{s}$ were defined in
(18).  We may then apply the above procedure to equation (34) and
obtain $\partial_{x}^{s}u\in H^{(0,2)}(\mu_{2}X)$,
$s=1,\ldots,m-2$.\par
   Lastly denote the right-hand side of (34) by $f_{s}$, then
the following equation holds in $L^{2}(\mu_{2}X)$,
\begin{equation}
z_{yy}=f_{s}-Az_{xx}-(D+sA_{x})z_{x}-Ez_{y}-(F+sD_{x}
+\frac{s(s-1)}{2}A_{xx})z.
\end{equation}
Since the right-hand side of (35) is in $H^{(0,1)}(\mu_{2}X)$, it
follows that $z_{yy}\in H^{(0,1)}(\mu_{2}X)$.  Then by
differentiating (35) with respect to $y$, we may apply a
boot-strap argument to obtain $u\in H^{m}(\mu_{2}X)$.  q.e.d.
\bigskip
\begin{center}
\textbf{3.  The Moser Estimate}
\end{center}
\bigskip

   Having established the existence of regular solutions to a
small perturbation of the linearized equation for (6), we intend
to apply a Nash-Moser type iteration procedure in the following
section, to obtain a smooth solution of (6) in $X$.  In the
current section, we shall make preparations for the Nash-Moser
procedure by establishing a certain a priori estimate.  This
estimate, referred to as the Moser estimate, will establish the
dependence of the solution $u$ of (22), on the coefficients of the
linearization as well as on the right-hand side, $f$.  If the
linearization is evaluated at some function $w\in
C^{\infty}(\mu_{2}\overline{X})$, then the Moser estimate is of
the form
\begin{equation}
\parallel u\parallel_{H^{m}}\leq C_{m}(\parallel f\parallel
_{H^{m}}+\parallel w\parallel_{H^{m+m_{1}}}\parallel
f\parallel_{H^{2}}),
\end{equation}
for some constants $C_{m}$ and $m_{1}$ independent of
$\varepsilon$ and $w$.\par
   Estimate (36) will first be established in the coordinates
$(\alpha,\beta)$, which we have been denoting by $(x,y)$ for
convenience, and later converted into the original coordinates
$(x,y)$ of the introduction.  We will need the Gagliardo-Nirenberg
estimates contained in the following lemma.
\medskip

\textbf{Lemma 3.1.}  \textit{Let $u,v\in
C^{k}(\overline{\Omega})$.}

$i)$\textit{  If $\sigma$ and $\varrho$ are multi-indices such
that $|\sigma|+|\varrho|=k$, then there exist constants
$\mathcal{M}_{1}$ and $\mathcal{M}_{2}$ depending on $k$, such
that}
\begin{equation*}
\parallel\partial^{\sigma}u\partial^{\varrho}v\parallel_{L^{2}(\Omega)}\leq
\mathcal{M}_{1}(|u|_{L^{\infty}(\Omega)}\parallel
v\parallel_{H^{k}(\Omega)}+\parallel
u\parallel_{H^{k}(\Omega)}|v|_{L^{\infty}(\Omega)}),
\end{equation*}
\textit{and}
\begin{equation*}
|\partial^{\sigma}u\partial^{\varrho}v|_{C^{0}(\overline{\Omega})}\leq
\mathcal{M}_{2}(|u|_{C^{0}(\overline{\Omega})}|v|_{C^{k}(\overline{\Omega})}+
|u|_{C^{k}(\overline{\Omega})}|v|_{C^{0}(\overline{\Omega})}).
\end{equation*}

$ii)$\textit{  If $\sigma_{1},\ldots,\sigma_{l}$ are multi-indices
such that $|\sigma_{1}|+\cdots+|\sigma_{l}|=k$, then there exists
a constant $\mathcal{M}_{3}$ depending on $l$ and $k$, such that}
\begin{eqnarray*}
&
&\!\!\parallel\partial^{\sigma_{1}}u_{1}\cdots\partial^{\sigma_{l}}u_{l}
\parallel_{L^{2}(\Omega)}\\
&\leq&\!\!
\mathcal{M}_{3}\sum_{j=1}^{l}(|u_{1}|_{L^{\infty}(\Omega)}\cdots
\widehat{|u_{j}|}_{L^{\infty}(\Omega)}
\cdots|u_{l}|_{L^{\infty}(\Omega)})\parallel
u_{j}\parallel_{H^{k}(\Omega)},
\end{eqnarray*}
\textit{where $\widehat{|u_{j}|}_{L^{\infty}(\Omega)}$ indicates
the absence of $|u_{j}|_{L^{\infty}(\Omega)}$.}

$iii)$\textit{  Let $B\subset\mathbb{R}^{N}$ be compact and
contain the origin, and let $G\in C^{\infty}(B)$.  If $u\in
H^{k+2}(\Omega,B)$ and $\parallel u\parallel _{H^{2}(\Omega)}\leq
\mathcal{C}$ for some fixed $\mathcal{C}$, then there exist
constants $\mathcal{M},\mathcal{M}_{k}>0$ such that}
\begin{equation*}
\parallel G\circ u\parallel_{H^{k}(\Omega)}\leq \mathcal{M}+\mathcal{M}_{k}\parallel
u\parallel_{H^{k+2}(\Omega)},
\end{equation*}
\textit{where }$\mathcal{M}=$Vol$(\Omega)|G(0)|$.\medskip

\textit{Proof.}  These estimates are standard consequences of the
interpolation inequalities, and may be found in, for instance,
[20]. q.e.d.\par\medskip
   Estimate (36) will follow by induction from the next two
propositions.  The first shall establish an estimate for the
$x$-derivatives only, while the second deals with all remaining
derivatives.\medskip

\textbf{Proposition 3.1.}  \textit{Suppose that the linearization,
$L_{1}$, is evaluated at some function $w\in
C^{\infty}(\mathbb{R}^{2})$ with $|w|_{C^{16}}\leq C_{1}$, as in}
(7).  \textit{Let $f\in H^{m}(\mathbb{R}^{2})$ and $u\in
H^{(m,1)}(\mathbb{R}^{2})\cap H^{m}(\mu_{2}X)$, $m\leq r-7$, be
the solution of} (22).  \textit{If $\varepsilon=\varepsilon(m)$ is
sufficiently small, then}
\begin{eqnarray*}
&
&\!\!\parallel\partial^{m}_{x}u\parallel+\parallel\partial^{m}_{x}u_{y}\parallel\\
&\leq&\!\! C_{m}(\parallel f\parallel_{m}+\parallel
u\parallel_{H^{m-1}(\mu_{2}X)}+\parallel
w\parallel_{H^{m+7}(\mu_{2}X)}\parallel
f\parallel_{H^{2}(\mu_{2}X)}),
\end{eqnarray*}
\textit{for some constant $C_{m}$ independent of $\varepsilon$ and
$w$.}\medskip

\textit{Proof.}   We proceed by induction on $m$.  The case $m=0$
is given by Lemma 2.1.  Now assume that the estimate holds for all
positive integers less than $m$.  Differentiate $L(w)u=f$
$m$-times with respect to $x$ and put $v=\partial_{x}^{m}u$, then
\begin{equation*}
L_{m}v=\partial_{x}^{m}f-\sum_{i=0}^{m-1}\partial_{x}^{i}(E_{x}
\partial_{x}^{m-1-i}u_{y}+\partial_{x}F_{m-1-i}\partial_{x}^{m-1-i}u)
:=f_{m},
\end{equation*}
where $L_{m}$ and $F_{m}$ were defined in (18).  If
$\varepsilon=\varepsilon(m)$ is sufficiently small, we can apply
Lemma 2.1 to obtain
\begin{equation}
\parallel\partial_{x}^{m}u\parallel+\parallel\partial_{x}^{m}u_{y}\parallel
\leq M\parallel f_{m}\parallel.
\end{equation}\par
   We now estimate each term of $f_{m}$. Let
$\parallel\cdot\parallel_{m,\text{ }\mu_{2}X}$ denote
$\parallel\cdot\parallel_{H^{m}(\mu_{2}X)}$, and
$|\cdot|_{\infty}$ denote $|\cdot|_{L^{\infty}(\mu_{2}X)}$.  A
calculation shows that
\begin{equation*}
\sum_{i=0}^{m-1}\partial_{x}^{i}(E_{x}\partial_{x}^{m-1-i}u_{y})=
mE_{x}\partial_{x}^{m-1}u_{y}+\sum_{i=1}^{m-1}\sum_{j=1}^{i}
\left(\begin{array}{c}i\\j\end{array}\right)\partial_{x}^{j+1}E
\partial_{x}^{m-1-j}u_{y}.
\end{equation*}
Then using Lemma 3.1 $(i)$ and $(iii)$, and recalling that $E_{x}$
vanishes on $\mathbb{R}^{2}-\mu_{2}X$, we obtain
\begin{eqnarray*}
&
&\parallel\sum_{i=0}^{m-1}\partial_{x}^{i}(E_{x}\partial_{x}^{m-1-i}u_{y})
\parallel\\
&\leq&\!\! M_{1}\parallel\partial_{x}^{m-1}u_{y}\parallel\\
&+&\!\! M_{2}(|\partial_{x}^{2}E|_{\infty}\parallel
u\parallel_{m-1,\text{
}\mu_{2}X}+\parallel\partial_{x}^{2}E\parallel_{m-1,\text{
}\mu_{2}X}|u|_{\infty})\\
&\leq&\!\! M_{1}\parallel\partial_{x}^{m-1}u_{y}\parallel\\
&+&\!\! M_{3}(|E|_{C^{2}(\mu_{2}\overline{X})}\parallel
u\parallel_{m-1,\text{ }\mu_{2}X}+\parallel w\parallel_{m+6,\text{
}\mu_{2}X}\parallel u
\parallel_{2,\text{ }\mu_{2}X}).
\end{eqnarray*}
Using the fact that $|E|_{C^{2}(\mu_{2}\overline{X})}\leq
C^{'}_{14}$ (Lemma 1.3), and the induction assumption, we have
\begin{equation}
\parallel\sum_{i=0}^{m-1}\partial_{x}^{i}(E_{x}\partial_{x}^{m-1-i}u_{y})
\parallel
\end{equation}
\begin{equation*}
\leq C^{'}_{m-1}(\parallel f\parallel_{m-1}+\parallel
u\parallel_{m-1,\text{ }\mu_{2}X}+\parallel w\parallel_{m+6,\text{
}\mu_{2}X}\parallel u\parallel_{2,\text{ }\mu_{2}X}).\nonumber
\end{equation*}\par
   In a similar manner, we may estimate
\begin{equation}
\parallel\sum_{i=0}^{m-1}\partial_{x}^{i}(\partial_{x}F_{m-1-i}\partial_{x}^{m-1-i}u)
\parallel
\end{equation}
\begin{equation*}
\leq C^{''}_{m-1}(\parallel f\parallel_{m-1}+\parallel
u\parallel_{m-1,\text{ }\mu_{2}X}+\parallel w\parallel_{m+7,\text{
}\mu_{2}X}\parallel u\parallel_{2,\text{ }\mu_{2}X}).
\end{equation*}
Furthermore, the methods used above can be made to show that
\begin{equation*}
\parallel u\parallel_{2,\text{ }\mu_{2}X}\leq
M_{4}\parallel f\parallel_{2,\text{ }\mu_{2}X}.
\end{equation*}
Then (38) and (39) yield
\begin{eqnarray*}
&
&\!\!\parallel\partial_{x}^{m}u\parallel+\parallel\partial_{x}^{m}u_{y}\parallel\\
&\leq&\!\! C_{m}(\parallel f\parallel_{m}+\parallel
u\parallel_{m-1,\text{ }\mu_{2}X}+\parallel w\parallel_{m+7,\text{
}\mu_{2}X}\parallel f\parallel_{2,\text{ }\mu_{2}X}),
\end{eqnarray*}
completing the proof by induction.  q.e.d.\par\medskip
   We now estimate the remaining derivatives.\medskip

\textbf{Proposition 3.2.}  \textit{Let $u$, $w$, $f$,
$\varepsilon$, and $m$ be as in Proposition} 3.1.  \textit{Then}
\begin{eqnarray*}
&
&\parallel\partial_{x}^{s}\partial_{y}^{t}u\parallel_{\mu_{2}X}\\
&\leq& C_{m}(\parallel f\parallel_{m,\text{ }\mu_{2}X}+\parallel
u\parallel_{m-1,\text{ }\mu_{2}X} +\parallel
w\parallel_{m+7,\text{ }\mu_{2}X}\parallel f\parallel_{2,\text{
}\mu_{2}X}),
\end{eqnarray*}
\textit{for }$0\leq s\leq m-t$, $0\leq t\leq m$, \textit{where
$C_{m}$ is independent of $\varepsilon$ and $w$.}\medskip

\textit{Proof.}  The cases $t=0,1$ are given by Proposition 3.1.
We will proceed by induction on $t$.  Assume that the desired
estimate holds for $0\leq s\leq m-t$, $0\leq t\leq k-1$, $0\leq
k\leq m$.\par
   Solving for $u_{yy}$ in the equation $L(w)u=f$, we obtain
\begin{equation}
u_{yy}=f-Au_{xx}-Du_{x}-Eu_{y}-Fu:=\overline{f}.
\end{equation}
Differentiate (40) with respect to
$\partial_{x}^{s}\partial_{y}^{k-2}$ where $0\leq s\leq m-k$, then
\begin{equation}
\partial_{x}^{s}\partial_{y}^{k}u=\partial_{x}^{s}\partial_{y}^{k-2}\overline{f}.
\end{equation}\par
   We now estimate each term on the right-hand side of (41).
Using Lemma 3.1 $(i)$ and $(iii)$, we have
\begin{eqnarray*}
& &\!\!\!\!\!\parallel\partial_{x}^{s}\partial_{y}^{k-2}(Au_{xx})\parallel_{\mu_{2}X}\\
&\leq&\!\!\!\!\!
M_{5}(\parallel\partial_{x}^{s+2}\partial_{y}^{k-2}
u\parallel_{\mu_{2}X}+\sum_{p\leq s,\text{ }q\leq
k-2\atop(p,q)\neq(0,0)}\parallel\partial_{x}^{p}\partial_{y}^{q}A
\partial_{x}^{s-p}\partial_{y}^{k-2-q}u_{xx}\parallel_{\mu_{2}X})\\
&\leq&\!\!\!\!\!
M^{'}_{5}(\parallel\partial_{x}^{s+2}\partial_{y}^{k-2}
u\parallel_{\mu_{2}X}\!+|A|_{C^{1}(\mu_{2}\overline{X})}\parallel\!
u\!\parallel _{m-1,\text{ }\mu_{2}X}\!+\parallel\!
A\!\parallel_{m,\text{
}\mu_{2}X}\!|u|_{\infty})\\
&\leq&\!\!\!\!\!
M^{''}_{5}(\parallel\partial_{x}^{s+2}\partial_{y}^{k-2}
u\parallel_{\mu_{2}X}+\parallel u\parallel_{m-1,\text{ }\mu_{2}X}+
\parallel w\parallel_{m+4,\text{ }\mu_{2}X}\parallel
f\parallel_{2,\mu_{2}X}).
\end{eqnarray*}
Furthermore, since $s\leq m-k$ the induction assumption implies
that
\begin{eqnarray*}
&
&\!\!\parallel\partial_{x}^{s+2}\partial_{y}^{k-2}u\parallel_{\mu_{2}X}\\
&\leq&\!\! M_{6}(\parallel f\parallel_{m,\text{ }\mu_{2}X}+
\parallel u\parallel_{m-1,\text{ }\mu_{2}X}+
\parallel w\parallel_{m+7,\text{ }\mu_{2}X}\parallel f\parallel
_{2,\text{ }\mu_{2}X}).
\end{eqnarray*}
Thus
\begin{eqnarray*}
&
&\!\!\parallel\partial_{x}^{s}\partial_{y}^{k-2}(Au_{xx})\parallel_{\mu_{2}X}\\
&\leq&\!\! M_{7}(\parallel f\parallel_{m,\text{ }\mu_{2}X}+
\parallel u\parallel_{m-1,\text{ }\mu_{2}X}+
\parallel w\parallel_{m+7,\text{ }\mu_{2}X}\parallel f\parallel
_{2,\text{ }\mu_{2}X}).
\end{eqnarray*}\par
   The remaining terms on the right-hand side of (41) may be
estimated in a similar manner.  Therefore
\begin{eqnarray*}
&
&\!\!\parallel\partial_{x}^{s}\partial_{y}^{k}u\parallel_{\mu_{2}X}\\
&\leq&\!\! M_{8}(\parallel f\parallel_{m,\text{ }\mu_{2}X}+
\parallel u\parallel_{m-1,\text{ }\mu_{2}X}+
\parallel w\parallel_{m+7,\text{ }\mu_{2}X}\parallel f\parallel
_{2,\text{ }\mu_{2}X}),
\end{eqnarray*}
for $0\leq s\leq m-k$.  The proof is now complete by induction.
q.e.d.\par\medskip
   By combining the previous two propositions,  we obtain the
following Moser estimate.\medskip

\textbf{Theorem 3.1.}  \textit{Let $u$, $w$, $f$, $\varepsilon$,
and $m$ be as in Proposition} 3.2.  \textit{Then}
\begin{equation*}
\parallel u\parallel_{m,\text{ }\mu_{2}X}\leq C_{m}(\parallel f\parallel
_{m,\text{ }\mu_{2}X} +\parallel w\parallel_{m+7,\text{
}\mu_{2}X}\parallel f\parallel_{2,\text{ }\mu_{2}X}),
\end{equation*}
\textit{where $C_{m}$ is independent of $\varepsilon$ and
$w$.}\medskip

\textit{Proof.}  This follows by induction on $m$, using
Proposition 3.2.  q.e.d.\par\medskip
   The Moser estimate of Theorem 3.1 is in terms of the variables
$(\alpha,\beta)$ of Lemma 1.3.  Since the Nash-Moser iteration
procedure of the following section will be carried out in the
original variables, $(x,y)$, of the introduction, we will now
obtain an analogous Moser estimate in these original coordinates.
Let $\parallel\cdot\parallel_{m,\text{ }\Omega}$,
$\parallel\cdot\parallel^{'}_{m,\text{ }\Omega}$, and
$\parallel\cdot\parallel^{''}_{m,\text{ }\Omega}$ denote the
$H^{m}(\Omega)$ norm with respect to the variables $(x,y)$,
$(\xi,\eta)$, and $(\alpha,\beta)$ respectively; a similar
notation will be used for the $C^{m}(\overline{\Omega})$ norms.
The following estimates will be needed in transforming the
estimate of Theorem 3.1 into the original variables.\medskip

\textbf{Lemma 3.2.}  \textit{If $\varepsilon=\varepsilon(m)$ is
sufficiently small, then}
\begin{equation*}
\parallel\xi_{x}\parallel_{m,\text{ }X_{1}}
\leq C_{m}(1+\parallel w\parallel_{m+7,\text{ }X_{1}}),
\end{equation*}
\textit{and}
\begin{equation*}
\parallel\alpha_{\xi}\parallel^{'}_{m,\text{ }X_{2}}
\leq C^{'}_{m}(1+\parallel w\parallel^{'}_{m+7,\text{ }X_{2}}),
\end{equation*}
\textit{where $C_{m}$ and $C_{m}^{'}$ are independent of
$\varepsilon$ and $w$, and $X_{1}$, $X_{2}$ were defined in
Lemmas} 1.2\textit{ and }1.3.\medskip

\textit{Proof.}  We shall only prove the first estimate, since a
similar argument yields the second.  The estimate will be proven
by induction on $m$. From the proof of Lemma 1.2 we have,
\begin{equation*}
|\xi_{x}|_{C^{0}(\overline{X}_{1})}\leq M_{9},
\end{equation*}
which gives the case $m=0$.  Now assume that the following
estimate holds,
\begin{equation*}
|\xi_{x}|_{C^{m-1}(\overline{X}_{1})}\leq C_{m-1}|b_{12}^{3}|
_{C^{m}(\overline{X}_{1})}.
\end{equation*}\par
   We will first estimate the $x$-derivatives.  Differentiate the
equation,
\begin{equation}
b_{12}^{3}(\xi_{x})_{x}+(\xi_{x})_{y}= -(b_{12}^{3})_{x}\xi_{x},
\end{equation}
$m$-times with respect to $x$ to obtain
\begin{equation*}
b_{12}^{3}(\partial_{x}^{m}\xi_{x})_{x}
+(\partial_{x}^{m}\xi_{x})_{y}=
-\partial_{x}^{m}[(b_{12}^{3})_{x}\xi_{x}]
-\sum_{i=0}^{m-1}\partial_{x}^{i}[(b_{12}^{3})_{x}
\partial_{x}^{m-i}\xi_{x}]:=g.
\end{equation*}
Then estimating $\partial_{x}^{m}\xi_{x}$ along the
characteristics of (42) as in the proof of Lemma 1.2, we find
\begin{equation*}
|\partial_{x}^{m}\xi_{x}|_{C^{0}(\overline{X}_{1})}
\leq\mu_{1}y_{0}|g|_{C^{0}(\overline{X}_{1})}.
\end{equation*}
Using the second half of Lemma 3.1 $(i)$ in the same way that the
first half was used in Proposition 3.1, and recalling that
$|b_{12}^{3}|_{C^{2}(\overline{X}_{1})}\leq\varepsilon M_{10}$,
produces
\begin{eqnarray*}
& &\!\!|g|_{C^{0}(\overline{X}_{1})}\\
&\leq&\!\! (m+1)\varepsilon M_{10}
|\partial_{x}^{m}\xi_{x}|_{C^{0}(\overline{X}_{1})}\\
& &\!\!+M_{10}^{'}(|(b_{12}^{3})_{xx}|
_{C^{0}(\overline{X}_{1})}|\xi_{x}|
_{C^{m-1}(\overline{X}_{1})}+|(b_{12}^{3})_{xx}|
_{C^{m-1}(\overline{X}_{1})}|\xi_{x}| _{C^{0}(\overline{X}_{1})}).
\end{eqnarray*}
Therefore if $\varepsilon$ is small enough to guarantee that
$(m+1)\mu_{1}y_{0}\varepsilon M_{10}<\frac{1}{2}$, we can bring
\begin{equation*}
(m+1)\mu_{1}y_{0}\varepsilon M_{10}
|\partial_{x}^{m}\xi_{x}|_{C^{0}(\overline{X}_{1})}
\end{equation*}
to the left-hand side:
\begin{equation}
|\partial_{x}^{m}\xi_{x}|_{C^{0}(\overline{X}_{1})}\leq
M_{11}(|\xi_{x}|_{C^{m-1}(\overline{X}_{1})}
+|b_{12}^{3}|_{C^{m+1}(\overline{X}_{1})}).
\end{equation}\par
   By solving for $(\xi_{x})_{y}$ in equation (42), and
differentiating the result with respect to
$\partial_{x}^{s}\partial_{y}^{t-1}$, $0\leq s\leq m-t$, $0\leq
t\leq m$, we can use the techniques of Proposition 3.2, combined
with Lemma 3.1 $(i)$ and (43), to obtain
\begin{equation}
|\partial_{x}^{s}\partial_{y}^{t}\xi_{x}|_{C^{0}(\overline{X}_{1})}
\leq M_{12}(|\xi_{x}|_{C^{m-1}(\overline{X}_{1})}
+|b_{12}^{3}|_{C^{m+1}(\overline{X}_{1})}).
\end{equation}
By the induction assumption on $m$, (44) implies that
\begin{equation*}
|\xi_{x}|_{C^{m}(\overline{X}_{1})}\leq
M_{13}|b_{12}^{3}|_{C^{m+1}(\overline{X}_{1})}.
\end{equation*}
Then the Sobolev Embedding Theorem gives
\begin{equation*}
\parallel\xi_{x}\parallel_{m,\text{ }X_{1}}\leq
M_{14}\parallel b_{12}^{3}\parallel_{m+3,\text{ }X_{1}}.
\end{equation*}
Thus, by Lemma 3.1 $(iii)$ we have
\begin{equation*}
\parallel\xi_{x}\parallel_{m,\text{ }X_{1}}\leq
M_{15}(1+\parallel w\parallel_{m+7,\text{ }X_{1}}).
\end{equation*}
q.e.d.\medskip

\textbf{Theorem 3.2.}  \textit{Let $u$, $w$, and $f$ be as in
Theorem} 3.1, \textit{and $m\leq r-25$. If
$\varepsilon=\varepsilon(m)$ is sufficiently small, then}
\begin{equation*}
\parallel u\parallel_{m,\text{ }X}\leq C_{m}(\parallel f\parallel
_{m,\text{ }X_{1}}+\parallel w\parallel _{m+25,\text{
}X_{1}}\parallel f\parallel _{2,\text{ }X_{1}}),
\end{equation*}
\textit{where $C_{m}$ is independent of $\varepsilon$ and
$w$.}\medskip

\textit{Proof.}  We first prove an analogue of the desired
estimate in terms of the variables $(\xi,\eta)$.  Observe that
\begin{equation}
\xi_{\alpha}=\frac{1}{\alpha_{\xi}}(\frac{\beta^{2}_{\eta}}
{\beta^{2}_{\eta}+\beta^{2}_{\xi}b_{12}^{5}})\geq M_{16}
\end{equation}
for some $M_{16}>0$, if $\varepsilon$ is sufficiently small.  Let
$G(b_{12}^{5})=\beta^{2}_{\eta}/(\beta^{2}_{\eta}+\beta^{2}_{\xi}b_{12}^{5})$,
and $s=m-t$, $0\leq t\leq m$.  A calculation shows that
\begin{equation*}
\parallel\partial_{\xi}^{s}\partial_{\eta}^{t}u\parallel^{'}
_{X_{2}}\leq M_{17}\sum_{k=0}^{m}\sum_{i=0}^{k}\parallel R_{ik}
\partial_{\alpha}^{k-i}\partial_{\beta}^{i}u\parallel_{\mu_{2}X}^{''},
\end{equation*}
where the $R_{ik}$ are polynomials in the variables
$\nabla_{\alpha,\beta}^{\sigma_{1}}\xi_{\alpha}$,
$\nabla_{\alpha,\beta}^{\sigma_{2}}\xi_{\alpha}^{-1}$,
$\nabla_{\alpha,\beta}^{\sigma_{3}}b_{12}^{5}$,
$\nabla_{\alpha,\beta}^{\sigma_{4}}G(b_{12}^{5})$,
$\nabla_{\xi,\eta}^{\sigma_{5}+1}\beta$, such that
$|\sigma_{j}|\leq m-k$, $1\leq j\leq 5$, and
$\sum_{\nu}|\sigma_{\nu}|\leq m-k$, where
$\sum_{\nu}|\sigma_{\nu}|$ represents the sum over all
$\sigma_{j}$ appearing in an arbitrary term of $R_{ik}$.  Then
using Lemma 3.1 $(ii)$ and $(iii)$, we find that
\begin{eqnarray}
\parallel\partial_{\xi}^{s}\partial_{\eta}^{t}u\parallel^{'}
_{X_{2}}\!\!\!\!&\leq&\!\!\!\! M_{18}[\parallel
u\parallel_{m,\text{ }\mu_{2}X}^{''}\nonumber\\
& &\!\!\!\!+(\parallel\xi_{\alpha}\parallel_{m,\text{
}\mu_{2}X}^{''} +\parallel\xi_{\alpha}^{-1}\parallel_{m,\text{
}\mu_{2}X}^{''} +\parallel b_{12}^{5}\parallel_{m+2,\text{
}\mu_{2}X}^{''})|u|_{\infty}]\\
&\leq&\!\!\!\! M_{18}^{'}[\parallel u\parallel_{m,\text{
}\mu_{2}X}^{''} +(\parallel\xi_{\alpha}\parallel_{m+2,\text{
}\mu_{2}X}^{''} +\parallel w\parallel_{m+6,\text{
}\mu_{2}X}^{''})|u|_{\infty}].\nonumber
\end{eqnarray}
Similarly,
\begin{equation}
\parallel\partial_{\alpha}^{s}\partial_{\beta}^{t}u\parallel^{''}
_{\mu_{2}X}\leq M_{19}[\parallel u\parallel^{'}_{m,\text{
}X_{2}}+(\parallel \alpha_{\xi}\parallel^{'}_{m+2,\text{
}X_{2}}+\parallel w\parallel^{'} _{m+6,\text{
}X_{2}})|u|_{\infty}].
\end{equation}
Then by Theorem 3.1, the Sobolev Lemma, and (46) we have
\begin{eqnarray}
\parallel\partial_{\xi}^{s}\partial_{\eta}^{t}u\parallel^{'}
_{X_{2}}\!\!&\leq&\!\! M_{20}(
\parallel f\parallel^{''}_{m,\text{ }\mu_{2}X}+\parallel
w\parallel^{''}_{m+7,\text{ }\mu_{2}X}\parallel f\parallel^{''}
_{2,\text{ }\mu_{2}X})\nonumber\\
& &\!\! +M^{'}_{20}(\parallel\xi_{\alpha}\parallel^{''}
_{m+2,\text{ }\mu_{2}X}+\parallel w\parallel^{''}_{m+6,\text{
}\mu_{2}X})\parallel f\parallel_{2,\text{ }X_{2}}^{'}.
\end{eqnarray}\par
   We now estimate the terms on the right-hand side of (48).
Using Lemma 3.1 $(i)$, $(iii)$, Lemma 3.2, (45), and (47) we have
\begin{eqnarray*}
& &\!\! \parallel\xi_{\alpha}\parallel^{''}_{m+2,\text{
}\mu_{2}X}\\
&\leq&\!\! M_{21}[\parallel\xi_{\alpha}\parallel_{m+2,\text{
}X_{2}}^{'}+(\parallel\alpha_{\xi}\parallel_{m+4,\text{
}X_{2}}^{'}+\parallel w\parallel_{m+8,\text{
}X_{2}}^{'})|\xi_{\alpha}|_{\infty}]\\
&\leq&\!\!
M_{22}[\parallel\alpha_{\xi}^{-1}G(b_{12}^{5})\parallel^{'}
_{m+2,\text{ }X_{2}}+\parallel\alpha_{\xi}\parallel_{m+4,\text{
}X_{2}}^{'}+\parallel w\parallel_{m+8,\text{
}X_{2}}^{'}]\\
&\leq&\!\!
M_{23}[|G(b_{12}^{5})|_{\infty}\parallel\alpha_{\xi}^{-1}\parallel^{'}
_{m+2,\text{ }X_{2}}+\parallel G(b_{12}^{5})\parallel
_{m+2,\text{ }X_{2}}^{'}|\alpha_{\xi}^{-1}|_{\infty}\\
& &\!\! +\parallel\alpha_{\xi}\parallel_{m+4,\text{
}X_{2}}^{'}+\parallel w\parallel_{m+8,\text{ }X_{2}}^{'}]\\
&\leq&\!\! M_{24}[\parallel\alpha_{\xi}\parallel_{m+4,\text{
}X_{2}}^{'}+\parallel w\parallel_{m+8,\text{ }X_{2}}^{'}]\\
&\leq&\!\! M_{25}[1+\parallel w\parallel^{'}_{m+11,\text{
}X_{2}}].
\end{eqnarray*}
Furthermore by (47), Lemma 3.2, and the Sobolev Lemma
\begin{eqnarray*}
& &\!\!\parallel f\parallel^{''}_{m,\text{ }\mu_{2}X}\\
&\leq&\!\! M_{26} [\parallel f\parallel^{'}_{m,\text{ }X_{2}}+
(\parallel\alpha_{\xi}\parallel^{'}_{m+2,\text{ }X_{2}}+\parallel
w\parallel^{'}_{m+6,\text{ }X_{2}})\parallel
f\parallel^{'}_{2,\text{ }X_{2}}]\\
&\leq&\!\! M_{26}^{'}[\parallel f\parallel^{'}_{m,\text{ }X_{2}}+
\parallel
w\parallel^{'}_{m+9,\text{ }X_{2}}
\parallel
f\parallel^{'}_{2,\text{ }X_{2}}].
\end{eqnarray*}
Also, the same method yields
\begin{eqnarray*}
\parallel w\parallel^{''}_{m+7,\text{ }\mu_{2}X}\!\!&\leq&\!\! M_{27}
(\parallel w\parallel^{'}_{m+7,\text{ }X_{2}}+
\parallel
w\parallel^{'}_{m+16,\text{ }X_{2}}
\parallel
w\parallel^{'}_{2,\text{ }X_{2}})\\
&\leq&\!\! M_{27}^{'}\parallel w\parallel^{'}_{m+16,\text{
}X_{2}}.
\end{eqnarray*}\par
   Therefore, from (48) and the above estimates we obtain
\begin{equation}
\parallel u\parallel^{'}_{m,\text{ }X_{2}}
\leq M_{28}(\parallel f\parallel^{'}_{m,\text{ }X_{2}}+\parallel
w\parallel^{'}_{m+16,\text{ }X_{2}}
\parallel f\parallel^{'}_{2,\text{
}X_{2}}).
\end{equation}
We can now apply the same procedure to obtain the following
analogue of (49) in terms of the original variables $(x,y)$,
\begin{equation*}
\parallel u\parallel_{m,\text{ }X}\leq M_{29}(\parallel f\parallel
_{m,\text{ }X_{1}}+\parallel w\parallel _{m+25,\text{
}X_{1}}\parallel f\parallel _{2,\text{ }X_{1}}).
\end{equation*}
q.e.d.
\bigskip
\begin{center}
\textbf{4.  The Nash-Moser Procedure}
\end{center}
\bigskip

   In this section we will carry out a Nash-Moser type iteration
procedure to obtain a solution of
\begin{equation}
\Phi(w)=0\text{ }\text{ in }\text{ }X.
\end{equation}
Instead of solving the linearized equation at each iteration, we
shall solve a small perturbation of the modified linearized
equation $L_{7}(v)u=f$, and then estimate the error at each step.
However, the theory of sections $\S 2$ and $\S 3$ requires that
$v$ and $f$ be defined on the whole plane.  Therefore, we will
need the following extension theorem.\medskip

\textbf{Theorem 4.1 [19].}  \textit{Let $\Omega$ be a bounded
convex domain in $\mathbb{R}^{2}$, with Lipschitz smooth boundary.
Then there exists a linear operator
$T_{\Omega}:L^{2}(\Omega)\rightarrow L^{2}(\mathbb{R}^{2})$ such
that:}

$i)$  $T_{\Omega}(g)|_{\Omega}=g$,

$ii)$  $T_{\Omega}:H^{m}(\Omega)\rightarrow H^{m}(\mathbb{R}^{2})$
\textit{continuously for each $m\in\mathbb{Z}_{\geq 0}$}.\medskip

   As with all Nash-Moser iteration schemes we will need smoothing
operators, which we now construct.  Fix $\widehat{\chi}\in
C^{\infty}_{c}(\mathbb{R}^{2})$ such that $\widehat{\chi}\equiv 1$
inside $X$.  Let
$\chi(x)=\int\int_{\mathbb{R}^{2}}\widehat{\chi}(\eta)e^{2\pi
i\eta\cdot x}d\eta$ be the inverse Fourier transform of
$\widehat{\chi}$.  Then $\chi$ is a Schwartz function and
satisfies $\int\int_{\mathbb{R}^{2}}\chi(x)dx\equiv 1$,
$\int\int_{\mathbb{R}^{2}}x^{\beta}\chi(x)dx=0$ for any
multi-index $\beta$, $\beta\neq 0$.  If $g\in
L^{2}(\mathbb{R}^{2})$ and $\mu\geq 1$, we define smoothing
operators $S^{'}_{\mu}:L^{2}(\mathbb{R}^{2})\rightarrow
H^{\infty}(\mathbb{R}^{2})$ by
\begin{equation*}
(S^{'}_{\mu}g)(x)=\mu^{2}\int\int
_{\mathbb{R}^{2}}\chi(\mu(x-y))g(y)dy.
\end{equation*}
Then we have (see [18]),\medskip

\textbf{Lemma 4.1.}  \textit{Let $l,m\in\mathbb{Z}_{\geq 0}$ and
$g\in H^{l}(\mathbb{R}^{2})$, then}

$i)$ $\parallel S_{\mu}^{'}g\parallel_{H^{m}(\mathbb{R}^{2})}\leq
C_{l,m}\parallel g\parallel_{H^{l}(\mathbb{R}^{2})}$, $m\leq l$,

$ii)$ $\parallel S_{\mu}^{'}g\parallel_{H^{m}(\mathbb{R}^{2})}\leq
C_{l,m}\mu^{m-l}\parallel g\parallel_{H^{l}(\mathbb{R}^{2})}$,
$l\leq m$,

$iii)$ $\parallel g-
S_{\mu}^{'}g\parallel_{H^{m}(\mathbb{R}^{2})}\leq
C_{l,m}\mu^{m-l}\parallel g\parallel_{H^{l}(\mathbb{R}^{2})}$,
$m\leq l$.\medskip

\noindent Furthermore, we obtain smoothing operators on $X$,
$S_{\mu}:L^{2}(X)\rightarrow H^{\infty}(X)$, by setting
$S_{\mu}g=(S^{'}_{\mu}Tg)|_{X}$, where $T$ is the extension
operator given by Theorem 4.1 with $\Omega=X$.  Moreover, it is
clear that the corresponding results of Lemma 4.1 hold for
$S_{\mu}$.\par
   We now set up the underlying iterative procedure.  Let
$\mu_{k}=\mu^{k}$, $S^{'}_{k}=S^{'}_{\mu_{k}}$,
$S_{k}=S_{\mu_{k}}$, and $w_{0}=0$.  Suppose that functions
$w_{0},w_{1},\ldots,w_{k}$ have been defined on $X$, and put
$v_{j}=S^{'}_{j}Tw_{j}$, $0\leq j\leq k$.  Let $L(v_{k})$ denote
the linearization of (50) evaluated at $v_{k}$, and let
$L_{8}(v_{k})$ be a small perturbation (on $X$) of $L_{7}(v_{k})$
to be given below, where $L_{7}(v_{k})$ is as in section $\S 1$.
Then define $w_{k+1}=w_{k}+u_{k}$ where $u_{k}$ is the solution,
restricted to $X$, of
\begin{equation}
L_{8}(v_{k})u_{k}=f_{k},
\end{equation}
given by Theorem 2.2 (see Lemma 4.2 below), and where $f_{k}$ will
be specified below.\par
   Let $Q_{k}(w_{k},u_{k})$ denote the quadratic error in the
Taylor expansion of $\Phi$ at $w_{k}$.  Then using the definition
of $L_{7}$ we have
\begin{eqnarray}
& &\!\!\!\Phi(w_{k+1})\\
\!\!\!&=&\!\!\!\Phi(w_{k})+L(w_{k})u_{k}+Q_{k}(w_{k},u_{k})\nonumber\\
&=&\!\!\!\Phi(w_{k})+A_{k}(w_{k})\partial_{xx}u_{k}+Q_{k}(w_{k},u_{k})\nonumber\\
&
&\!\!\!+\varepsilon(1\!+\varepsilon(w_{k})_{xx}\!+\varepsilon^{2n}
H^{n}P_{11}(w_{k}))(P^{6}_{22}(w_{k})L_{7}(w_{k})u_{k}\!+D_{k}(w_{k})
\partial_{x}u_{k})\nonumber\\
&=&\!\!\!\Phi(w_{k})\!+\!\varepsilon(1+\varepsilon(v_{k})_{xx}+\varepsilon^{2n}
H^{n}P_{11}(v_{k}))|_{X}P^{6}_{22}(v_{k}|_{X})L_{8}(v_{k}|_{X})u_{k}\nonumber\\
& &\!\!\!+e_{k},\nonumber
\end{eqnarray}
where
\begin{eqnarray*}
e_{k}\!\!\!&=&\!\!\!\varepsilon(P_{k}(w_{k})L_{8}(w_{k})-P_{k}(v_{k}|_{X})L_{8}(v_{k}|_{X}))u_{k}
+A_{k}(w_{k})\partial_{xx}u_{k}\\
&
&\!\!\!+Q_{k}(w_{k},u_{k})-\varepsilon\overline{P}_{k}(w_{k})(P^{6}_{22}(w_{k})\overline{A}_{k}
\partial_{\alpha\alpha}u_{k}-(S_{k}D_{k}(w_{k}))\partial_{x}u_{k}),
\end{eqnarray*}
\begin{equation*}
P_{k}(w_{k})=(1+\varepsilon(w_{k})_{xx}+\varepsilon^{2n}H^{n}P_{11}(w_{k}))
P^{6}_{22}(w_{k}),
\end{equation*}
\begin{equation*}
\overline{P}_{k}(w_{k})=1+\varepsilon(w_{k})_{xx}+\varepsilon^{2n}H^{n}P_{11}(w_{k}),
\end{equation*}
\begin{equation*}
A_{k}(w_{k})=\varepsilon\overline{P}_{k}^{-1}(w_{k})\Phi(w_{k}),\text{
}\text{ }\text{ }\text{ }\text{
}\overline{A}_{k}=\varepsilon^{n}\mu_{k}^{-4}\beta\phi(\alpha)\phi(\beta)+\psi_{1}(\beta),
\end{equation*}
\begin{equation*}
D_{k}(w_{k})=\frac{1}{2}\partial_{x}[\overline{P}_{k}^{-2}(w_{k})\Phi(w_{k})]
+\frac{1}{2}\overline{P}_{k}^{-2}(w_{k})\partial_{x}\Phi(w_{k}),
\end{equation*}
\begin{eqnarray*}
L_{8}(w_{k})u_{k}\!\!\!&=&\!\!\!L_{7}(w_{k})u_{k}+\overline{A}_{k}\partial_{\alpha\alpha}u_{k}\\
&
&\!\!\!+\phi(\alpha)\phi(\beta)T[(P^{6}_{22}(w_{k}))^{-1}(I-S_{k})D_{k}(w_{k})]\partial_{x}u_{k},
\end{eqnarray*}
the functions $\phi$ and $\psi_{1}$ are as in section $\S 2$,
$(\alpha,\beta)$ are the coordinates of Lemma 1.3; note also also
that we use $\phi|_{X}\equiv 1$ and $T(\cdot)|_{X}=I$ in (52).\par
   We now define $f_{k}$.  In order to solve (51) with the theory of
section $\S 2$, we require $f_{k}$ to be defined on all of
$\mathbb{R}^{2}$.  Furthermore, we need the right-hand side of
(52) to tend to zero sufficiently fast, to make up for the error
incurred at each step by solving (51) instead of solving the
unmodified linearized equation.  Therefore we set $E_{0}=0$,
$E_{k}=\sum_{i=0}^{k-1}e_{i}$, and define
\begin{equation*}
f_{0}=-T[(\varepsilon P_{0}(v_{0}))^{-1}S_{0}\Phi(w_{0})],
\end{equation*}
\begin{equation*}
f_{k}=T[(\varepsilon
P_{k}(v_{k}))^{-1}(S_{k-1}E_{k-1}-S_{k}E_{k}+(S_{k-1}-S_{k})\Phi(w_{0}))].
\end{equation*}
It follows that
\begin{eqnarray}
\Phi(w_{k+1})&=&\Phi(w_{0})+\sum_{i=0}^{k}
\varepsilon P_{i}(v_{i}|_{X})(f_{i}|_{X})+E_{k}+e_{k}\\
&=&(I-S_{k})\Phi(w_{0})+(I-S_{k})E_{k}+e_{k}.\nonumber
\end{eqnarray}
In what follows, we will show that the right-hand side of (53)
tends to zero sufficiently fast to guarantee the convergence of
$\{w_{k}\}_{k=0}^{\infty}$ to a solution of (50).\par
   Let $b$ be a positive number that will be chosen as large as
possible, set $\delta=\varepsilon^{n-1}$, and
$\mu=\varepsilon^{\frac{1-n}{b+1}}$.  Furthermore, let $m_{*}\in
\mathbb{Z}_{\geq 0}$ be such that $\Phi(w_{0})\in H^{m_{*}}(X)$.
For convenience we will denote the $H^{m}(X)$ and
$H^{m}(\mathbb{R}^{2})$ norms by $\parallel\cdot\parallel_{m}$ and
$\parallel\cdot\parallel_{m,\text{ }\mathbb{R}^{2}}$,
respectively.  The convergence of $\{w_{k}\}_{k=0}^{\infty}$ will
follow from the following eight statements, valid for $0\leq m\leq
m_{*}-25$ unless specified otherwise, which shall be proven by
induction on $j$, for some constants $C_{1},C_{2},C_{3}$, and
$C_{4}$ independent of $j$, $\varepsilon$, and $\mu$, but
dependent on $m$.\bigskip

   I$_{j}$:   $\parallel u_{j-1}\parallel_{m}\leq\delta\mu_{j-1}^{m-b}$,\bigskip

   II$_{j}$:   $\parallel w_{j}\parallel_{m}\leq\begin{cases}
C_{1}\delta & \text{if $m-b\leq-1/2$},\\
C_{1}\delta\mu_{j}^{m-b} & \text{if $m-b\geq 1/2$,}
\end{cases}$\bigskip

   III$_{j}$:   $\parallel w_{j}\parallel_{18}\leq C_{1}\delta,\text{ }
   \text{ }\parallel v_{j}\parallel_{18,\text{ }\mathbb{R}^{2}}\leq
   C_{3}\delta$,\bigskip

   IV$_{j}$:   $\parallel w_{j}-v_{j}\parallel_{m}\leq C_{2}\delta
   \mu_{j}^{m-b}$,\bigskip

   V$_{j}$:   $\parallel v_{j}\parallel_{m,\text{ }\mathbb{R}^{2}}\leq\begin{cases}
C_{3}\delta & \text{if $m-b\leq-1/2$},\\
C_{3}\delta\mu_{j}^{m-b} & \text{if $m-b\geq 1/2$,}
\end{cases}$ $0\leq m<\infty$,\bigskip

   VI$_{j}$:   $\parallel e_{j-1}\parallel_{m}\leq
   \varepsilon\delta^{2}\mu_{j-1}^{m-b}$, $0\leq m\leq m_{*}-30$,\bigskip

   VII$_{j}$:   $\parallel f_{j}\parallel_{m,\text{ }\mathbb{R}^{2}}
   \leq C_{4}\delta^{2}(1+\mu^{b-m})\mu_{j}^{m-b}$, $0\leq m\leq m_{*}$,\bigskip

   VIII$_{j}$:   $\parallel\Phi(w_{j})\parallel_{m}
   \leq\delta\mu_{j}^{m-b}$, $0\leq m\leq m_{*}-30$.\bigskip

   Assume that the above eight statements hold for $j=0,\ldots,k$.
Before showing the induction step we will need the following
preliminary lemma which allows us to study equation (51).\medskip

\textbf{Lemma 4.2.}  \textit{If $\varepsilon$ is sufficiently
small, then the theory of sections $\S 2$ and $\S 3$ applies to
the operators $L_{8}(v_{k})$ and $L_{8}(v_{0})$.}
\medskip

\textit{Proof.}  We first show that Lemma 2.1 holds for
$L_{8}(v_{k})$.  Extend the coefficients of $L_{7}(v_{k})$ to the
entire $\alpha\beta$-plane and denote them by $A_{k}$, $D_{k}$,
$E_{k}$, $F_{k}$ as in section $\S 2$.  Write
\begin{equation*}
L_{8}(v_{k})=\widetilde{A}_{k}\partial_{\alpha\alpha}+\partial_{\beta\beta}
+\widetilde{D}_{k}\partial_{\alpha}+\widetilde{E}_{k}\partial_{\beta}
+\widetilde{F}_{k},
\end{equation*}
let $I_{i}$, $i=1,2,3,4$, be as in the proof of Lemma 2.1, and let
$\widetilde{I}_{i}$ be analogous to $I_{i}$ with $A_{k}$, $D_{k}$,
$E_{k}$, $F_{k}$ replaced by $\widetilde{A}_{k}$,
$\widetilde{D}_{k}$, $\widetilde{E}_{k}$, $\widetilde{F}_{k}$.
Then a calculation shows that
\begin{equation*}
\widetilde{I}_{1}\geq I_{1}+
\begin{cases}
\varepsilon\delta\mu_{k}^{-4}\phi(\alpha)(\frac{1}{2}+O(|\beta|))
& \text{if }|\beta|\leq y_{3},\\
C+O(\varepsilon) & \text{if }|\beta|\geq y_{3},
\end{cases}
\end{equation*}
for some constant $C>0$ independent of $\varepsilon$ and $k$,
where $y_{3}$ is as in the proof of Lemma 2.1.  Furthermore, using
the definition of $\Phi$, Lemma 3.1 $(iii)$, and III$_{k}$, we
have
\begin{eqnarray*}
|(I-S_{k})D_{k}(v_{k})|_{C^{0}(X)}&\leq&
C\parallel(I-S_{k})D_{k}(v_{k})\parallel_{2}\\
&\leq& C\mu_{k}^{-5}\parallel D_{k}(v_{k})\parallel_{7}\\
&\leq& C\mu_{k}^{-5}(\varepsilon\parallel
v_{k}\parallel_{12}+\varepsilon^{2n})\\
&\leq& C\varepsilon\delta\mu_{k}^{-5}
\end{eqnarray*}
since $\Phi(0)=O(\varepsilon^{2n})$.  It follows that
\begin{equation*}
\widetilde{I}_{3}\geq
I_{3}+O(\varepsilon\delta\mu_{k}^{-5}\phi(\alpha)),\text{ }\text{
}\text{ }\text{ }\text{ }\text{ }\widetilde{I}_{4}=
I_{4}+O(\varepsilon),
\end{equation*}
\begin{equation*}
\widetilde{I}_{2}=I_{2}+O(\varepsilon\delta\mu_{k}^{-4}|\phi^{'}(\alpha)|
+\varepsilon\delta\mu_{k}^{-5}\phi(\alpha)),
\end{equation*}
from which we also find
\begin{equation*}
\widetilde{I}_{1}\widetilde{I}_{3}-2\widetilde{I}_{2}^{2}\geq
I_{1}I_{3}-2I_{2}^{2}+\varepsilon\delta\mu_{k}^{-4}\phi(\alpha)
(C+O(\mu_{k}^{-1}+\varepsilon))\geq 0,
\end{equation*}
if $\varepsilon$ is sufficiently small.  We then conclude that
Lemma 2.1 holds for $L_{8}(v_{k})$.  Similarly, the proofs of the
remaining results of sections $\S 2$ and $\S 3$ need only slight
modifications to show that they also hold for $L_{8}(v_{k})$.
Lastly, the same method applies to $L_{8}(v_{0})$ if we note that
\begin{equation*}
|(I-S_{0})D_{0}(v_{0})|_{C^{0}(X)}\leq C\varepsilon^{2n}.
\end{equation*}
q.e.d.\par\medskip
   The next four propositions will show that the above eight
statements hold for $j=k+1$.  The case $j=0$ will be proven
shortly there after.\medskip

\textbf{Proposition 4.1.}  \textit{If $27\leq b\leq m_{*}-26$,
$0\leq m\leq m_{*}-25$, and $\varepsilon$ is sufficiently small,
then} I$_{k+1}$, II$_{k+1}$, III$_{k+1}$, IV$_{k+1}$, \textit{and}
V$_{k+1}$ \textit{hold.}
\medskip

\textit{Proof.}  I$_{k+1}$:  First note that by III$_{k}$,
\begin{equation*}
|v_{k}|_{C^{16}(\mathbb{R}^{2})}\leq C\parallel v_{k}
\parallel_{18,\text{ }\mathbb{R}^{2}}\leq C^{'}.
\end{equation*}
Therefore, we may apply Lemma 4.2 and the theory of section $\S 2$
to obtain the solution $u_{k}$ of (51).  We require $m\leq
m_{*}-25$ so that the hypotheses of Theorem 3.2 are fulfilled.  If
$m+25-b\geq1/2$ then using Theorem 3.2, V$_{k}$, VII$_{k}$, and
$b\geq 27$, we have
\begin{eqnarray*}
\parallel u_{k}\parallel_{m}\!\!\!&\leq&\!\!\! C_{m}(\parallel
f_{k}\parallel_{m,\text{ }\mathbb{R}^{2}}+\parallel
v_{k}\parallel_{m+25,\text{ }\mathbb{R}^{2}}\parallel
f_{k}\parallel_{2,\text{ }\mathbb{R}^{2}})\\
&\leq&\!\!\!
C_{m}(C_{4}\delta^{2}(1+\mu^{b-m})\mu_{k}^{m-b}+C_{3}C_{4}\delta^{3}
(1+\mu^{b-2})\mu_{k}^{m+25-b}\mu_{k}^{2-b})\\
&\leq&\!\!\!\delta\mu_{k}^{m-b},
\end{eqnarray*}
if $\varepsilon$ is sufficiently small, since
$\delta\mu^{b-m}=\varepsilon^{(n-1)(1-\frac{b-m}{b+1})}\leq\varepsilon^{\frac{1}{b+1}}$.
If $m+25-b\leq-1/2$ and $m\geq2$, then using $\parallel
v_{k}\parallel_{m+25,\text{ }\mathbb{R}^{2}}\leq C_{3}\delta$ in
the estimate above gives the desired result.  Furthermore, if
$0\leq m<2$ then the methods of Theorem 3.2 show that $\parallel
u_{k}\parallel_{m}\leq M\parallel f_{k}\parallel_{m,\text{
}\mathbb{R}^{2}}$; in which case VII$_{k}$ gives the desired
result.\par
   II$_{k+1}$:  Since $w_{k+1}=\sum_{i=0}^{k}u_{i}$, we have
\begin{equation*}
\parallel w_{k+1}\parallel_{m}\leq\sum_{i=0}^{k}\parallel u_{i}
\parallel_{m}\leq\delta\sum_{i=0}^{k}\mu_{i}^{m-b}.
\end{equation*}
Hence, if $m-b\leq-1/2$
\begin{equation*}
\parallel
w_{k+1}\parallel_{m}\leq\delta\sum_{i=0}^{\infty}(\mu^{i})^{-1/2}
\leq\delta\sum_{i=0}^{\infty}(2^{i})^{-1/2}:=C_{1}\delta,
\end{equation*}
and if $m-b\geq 1/2$,
\begin{equation*}
\parallel w_{k+1}\parallel_{m}\leq\delta\mu_{k+1}^{m-b}
\sum_{i=0}^{k}(\frac{\mu_{i}}{\mu_{k+1}})^{m-b}
\leq\delta\mu_{k+1}^{m-b}\sum_{i=0}^{\infty}(\mu^{-i})^{1/2} \leq
C_{1}\delta\mu_{k+1}^{m-b}.
\end{equation*}\par
   III$_{k+1}$:  Since $b\geq 27$ we have $18-b\leq-1/2$.
Therefore II$_{k+1}$ and V$_{k+1}$ (proven below) imply that
\begin{equation*}
\parallel w_{k+1}\parallel_{18}\leq C_{1}\delta
\text{ }\text{ }\text{ and }\text{ }
\parallel v_{k+1}\parallel_{18,\text{ }\mathbb{R}^{2}}\leq
C_{3}\delta.
\end{equation*}\par
   IV$_{k+1}$:  Since $b\leq m_{*}-26$ we have $m_{*}-25-b\geq 1/2$.
Therefore Lemma 4.1 and II$_{k+1}$ yield,
\begin{eqnarray*}
\parallel w_{k+1}-v_{k+1}\parallel_{m}&=&\parallel
(I-S_{k+1})w_{k+1}\parallel_{m}\\
&\leq& C_{m}\mu^{m-(m_{*}-25)}_{k+1}\parallel w_{k+1}\parallel_{m_{*}-25}\\
&\leq& C_{m}\mu^{m-(m_{*}-25)}_{k+1}C_{1}\delta\mu^{m_{*}-25-b}_{k+1}\\
&:=& C_{2}\delta\mu^{m-b}_{k+1}.
\end{eqnarray*}\par
   V$_{k+1}$:  From Lemma 4.1 and $b\leq m_{*}-26$ we have for all
$m\geq 0$,
\begin{eqnarray*}
\parallel v_{k+1}\parallel_{m,\text{ }\mathbb{R}^{2}}\!\!\!&=&\!\!\!\parallel
S^{'}_{k+1}Tw_{k+1}\parallel_{m,\text{ }\mathbb{R}^{2}}\\
&\leq&\!\!\! C^{'}_{m}\!\parallel T\parallel\!
\begin{cases}
\parallel w_{k+1}\parallel_{b-1}
& \text{if }m-b\leq-1/2,\\
\mu_{k+1}^{m-b-1}\!\parallel w_{k+1}\parallel_{b+1} & \text{if
}m-b\geq 1/2.
\end{cases}
\end{eqnarray*}
V$_{k+1}$ now follows from II$_{k+1}$.  q.e.d.\par\medskip
   Write $e_{k}=e^{'}_{k}+e^{''}_{k}+e^{'''}_{k}$, where
\begin{eqnarray*}
e^{'}_{k}\!\!\!&=&\!\!\!\varepsilon(P_{k}(w_{k})L_{8}(w_{k})-P_{k}(v_{k}|_{X})L_{8}(v_{k}|_{X}))u_{k},\\
e^{''}_{k}\!\!\!&=&\!\!\!-\varepsilon\overline{P}_{k}(w_{k})
(P^{6}_{22}(w_{k})\overline{A}_{k}\partial_{\alpha\alpha}u_{k}
-(S_{k}D_{k}(w_{k}))\partial_{x}u_{k})+A_{k}(w_{k})\partial_{xx}u_{k},\\
e^{'''}_{k}\!\!\!&=&\!\!\!Q_{k}(w_{k},u_{k}).
\end{eqnarray*}

\textbf{Proposition 4.2.}  \textit{If the hypotheses of
Proposition} 4.1 \textit{hold and $0\leq m\leq m_{*}-30$, then}
VI$_{k+1}$ \textit{holds.}\medskip

\textit{Proof.}  We will estimate $e^{'}_{k}$, $e^{''}_{k}$, and
$e^{'''}_{k}$ separately.  Denote
\begin{eqnarray*}
&
&\!\!\!(P_{k}(w_{k})L_{8}(w_{k})-P_{k}(v_{k}|_{X})L_{8}(v_{k}|_{X}))u_{k}\\
&=&\!\!\!\sum_{i,j}d_{ij}(u_{k})_{x_{i}x_{j}}
+\sum_{i}d_{i}(u_{k})_{x_{i}}+du_{k},
\end{eqnarray*}
then Lemma 3.1 $(i)$ and $(iii)$, I$_{k}$, and IV$_{k}$ show that
\begin{eqnarray*}
\parallel e_{k}^{'}\parallel_{m}&\leq& \varepsilon C_{m,1}
[(\sum_{i,j}\parallel d_{ij}\parallel_{m}+\sum_{i}\parallel
d_{i}\parallel_{m}+\parallel d\parallel_{m})\parallel
u_{k}\parallel_{4}\\
& &+(\sum_{i,j}\parallel d_{ij}\parallel_{2}+\sum_{i}\parallel
d_{i}\parallel_{2}+\parallel
d\parallel_{2})\parallel u_{k}\parallel_{m+2}]\\
&\leq& \varepsilon C_{m,2}(\parallel
w_{k}-v_{k}\parallel_{m+5}\parallel u_{k}\parallel_{4}+\parallel
w_{k}-v_{k}\parallel_{7}\parallel u_{k}\parallel_{m+2})\\
&\leq& C_{m,3}\varepsilon\delta^{2}\mu_{k}^{9-b}\mu_{k}^{m-b}\\
&\leq& \frac{\varepsilon}{3}\delta^{2}\mu_{k}^{m-b}
\end{eqnarray*}
if $\varepsilon$ is sufficiently small, since
$\mu_{k}^{9-b}\leq\mu^{9-b}=\varepsilon^{(9-b)(\frac{1-n}{b+1})}\leq
\varepsilon^{18/28}$.  Note that we have also used $m\leq
m_{*}-30$, which allows us to apply IV$_{k}$.\par
   We now estimate $e^{''}_{k}$.  By Lemma 3.1 $(i)$ and $(iii)$,
I$_{k}$, II$_{k}$, and VIII$_{k}$,
\begin{eqnarray*}
& &\!\!\!\parallel A_{k}\partial_{xx}u_{k}\parallel_{m}\\
&\leq&\!\!\! C_{m,4}(\parallel\partial_{xx}u_{k}\parallel_{2}
\parallel A_{k}\parallel_{m}+\parallel\partial_{xx}u_{k}\parallel_{m}
\parallel A_{k}\parallel_{2})\\
&\leq&\!\!\!\varepsilon C_{m,5}[\parallel
u_{k}\parallel_{4}((1+\parallel
w_{k}\parallel_{6})\parallel\Phi(w_{k})\parallel_{m}+\parallel
w_{k}\parallel_{m+4}\parallel\Phi(w_{k})\parallel_{2})\\
& &\!\!\!+\parallel u_{k}\parallel_{m+2}
\parallel\Phi(w_{k})\parallel_{2}]\\
&\leq&\!\!\!\varepsilon
C_{m,6}[\delta\mu^{4-b}_{k}(\delta\mu^{m-b}_{k}+\delta^{2}
\mu^{m+4-b}_{k}\mu_{k}^{2-b})+\delta^{2}\mu^{m+2-b}_{k}\mu^{2-b}_{k}]\\
&\leq&\!\!\!\varepsilon C_{m,7}\mu^{10-b}_{k}\delta^{2}\mu^{m-b}_{k}\\
&\leq&\!\!\!\frac{\varepsilon}{9}\delta^{2}\mu^{m-b}_{k},
\end{eqnarray*}
if $\varepsilon$ is sufficiently small and $m+4-b\geq 1/2$.  If
$m+4-b\leq -1/2$ then we may use the estimate $\parallel
w_{k}\parallel_{m+4}\leq C_{1}\delta$ to obtain the same outcome.
Furthermore, the same methods combined with Lemma 4.1 show that
\begin{eqnarray*}
&
&\!\!\!\parallel\!\varepsilon\overline{P}_{k}(w_{k})(S_{k}D_{k})\partial_{x}u_{k}\!\parallel_{m}\\
&\leq& \!\!\!\varepsilon
C_{m,8}(\parallel\!\partial_{x}u_{k}\!\parallel_{2}\parallel\!\overline{P}_{k}(S_{k}D_{k})
\!\parallel_{m}+\parallel\!\partial_{x}u_{k}\!\parallel_{m}
\parallel\!\overline{P}_{k}(S_{k}D_{k})\!
\parallel_{2})\\
&\leq& \!\!\!\varepsilon C_{m,9}[\parallel\!
u_{k}\!\parallel_{3}(\mu_{k}
\parallel\!\overline{P}_{k}\!\parallel_{2}\parallel \!D_{k}\!\parallel_{m-1}+
\parallel\!\overline{P}_{k}\!\parallel_{m}\parallel \!D_{k}\!\parallel_{2})\\
& &\!\!\!+\parallel \!u_{k}\!\parallel_{m+1}\parallel \!D_{k}
\!\parallel_{2}]\\
&\leq& \!\!\!\varepsilon C_{m,10}[\parallel\!
u_{k}\!\parallel_{3}(\mu_{k}\!
\parallel\!\Phi(w_{k})\!\parallel_{m}\!+
\mu_{k}(1+\!\parallel \!w_{k}\!\parallel_{m+4})\!\parallel\!\Phi(w_{k})\!\parallel_{3})\\
& &\!\!\!+\parallel \!u_{k}\!\parallel_{m+1}\parallel\!\Phi(w_{k})\!\parallel_{3}]\\
&\leq& \!\!\!\varepsilon
C_{m,11}[\delta\mu_{k}^{3-b}(\delta\mu_{k}^{m+1-b}
+\delta^{2}\mu_{k}^{m+5-b}\mu_{k}^{3-b})+\delta^{2}\mu_{k}^{m+1-b}
\mu_{k}^{3-b}]\\
&\leq& \!\!\!\varepsilon C_{m,12}\mu_{k}^{11-b}\delta^{2}\mu_{k}^{m-b}\\
&\leq& \!\!\!\frac{\varepsilon}{9}\delta^{2}\mu^{m-b}_{k}.
\end{eqnarray*}
Similarly, since $\psi_{1}(\beta)\equiv 0$ in $X$ it follows that
\begin{eqnarray*}
& &\!\!\!\parallel\varepsilon
P_{k}(w_{k})\overline{A}_{k}\partial_{\alpha\alpha}u_{k}\parallel_{m}\\
&\leq&\!\!\! \varepsilon^{2}\delta\mu_{k}^{-4} C_{m,13}(\parallel
u_{k}\parallel_{4}\parallel w_{k}\parallel_{m+4}+\parallel
u_{k}\parallel_{m+2}(1+\parallel w_{k}\parallel_{6}))\\
&\leq&\!\!\! \varepsilon^{2}\delta\mu_{k}^{-4}
C_{m,14}(\delta^{2}\mu_{k}^{4-b}\mu_{k}^{m+4-b}+\delta\mu_{k}^{m+2-b})\\
&\leq&\!\!\! \frac{\varepsilon}{9}\delta^{2}\mu_{k}^{m-b}.
\end{eqnarray*}
Therefore
\begin{equation*}
\parallel e^{''}_{k}\parallel\leq
\frac{\varepsilon}{3}\delta^{2}\mu^{m-b}_{k}.
\end{equation*}\par
   We now estimate $e^{'''}_{k}$.  We have
\begin{equation*}
e^{'''}_{k}=Q_{k}(w_{k},u_{k})=\int_{0}^{1}(1-t)\frac{\partial^{2}}
{\partial t^{2}}\Phi(w_{k}+tu_{k})dt.
\end{equation*}
Apply Lemma 3.1 $(i)$ and $(ii)$, as well as the Sobolev Lemma to
obtain
\begin{eqnarray*}
\parallel e^{'''}_{k}\parallel_{m}&\leq&\int_{0}^{1}
\sum_{|\sigma|,|\gamma|\leq
2}\parallel\nabla_{\overline{\sigma}\overline{\gamma}}
\Phi(w_{k}+tu_{k})\partial^{\sigma}u_{k}\partial^{\gamma}u_{k}\parallel
_{m}dt\\
&\leq&\int_{0}^{1}\sum_{|\sigma|,|\gamma|\leq
2}C_{m,15}(|\nabla_{\overline{\sigma}\overline{\gamma}}\Phi(w_{k}+tu_{k})|_{\infty}
\parallel\partial^{\sigma}u_{k}\partial^{\gamma}u_{k}\parallel_{m}\\
&
&+\parallel\nabla_{\overline{\sigma}\overline{\gamma}}\Phi(w_{k}+tu_{k})\parallel_{m}
|\partial^{\sigma}u_{k}\partial^{\gamma}u_{k}|_{\infty})dt\\
&\leq&\int_{0}^{1}C_{m,16}(\parallel\nabla^{2}\Phi(w_{k}+tu_{k})\parallel_{2}
\parallel u_{k}\parallel_{4}\parallel u_{k}\parallel_{m+2}\\
& &+\parallel\nabla^{2}\Phi(w_{k}+tu_{k})\parallel_{m}
\parallel u_{k}\parallel_{4}^{2})dt,
\end{eqnarray*}
where $\overline{\sigma}=\partial^{\sigma}(w_{k}+tu_{k})$ and
$\overline{\gamma}=\partial^{\gamma}(w_{k}+tu_{k})$.  The notation
$\nabla^{2}\Phi$ represents the collection of second partial
derivatives with respect to the variables $\overline{\sigma}$,
$\overline{\gamma}$, so by (6)
$\nabla^{2}\Phi=O(\varepsilon^{2})$. Therefore using Lemma 3.1
$(iii)$, I$_{k}$, and II$_{k}$, we have
\begin{eqnarray*}
\parallel e^{'''}_{k}\parallel_{m}&\leq& \varepsilon^{2} C_{m,17}[(1+\parallel
w_{k}\parallel_{6}+\parallel u_{k}\parallel_{6})\parallel u_{k}
\parallel_{4}\parallel u_{k}\parallel_{m+2}\\
& &+(1+\parallel w_{k}
\parallel_{m+4}+\parallel u_{k}\parallel_{m+4})\parallel u_{k}
\parallel^{2}_{4}]\\
&\leq& \varepsilon^{2}
C_{m,18}[\delta^{2}\mu^{4-b}_{k}\mu^{m+2-b}_{k}+\delta^{2}\mu_{k}^{2(4-b)}
+\delta^{3}
\mu^{m+4-b}_{k}\mu^{2(4-b)}_{k}]\\
&\leq&\frac{\varepsilon}{3}\delta^{2}\mu^{m-b}_{k}
\end{eqnarray*}
if $\varepsilon$ is sufficiently small, since $b\geq 27$.
Combining the estimates of $e^{'}_{k}$, $e^{''}_{k}$, and
$e^{'''}_{k}$ yields the desired result.  q.e.d.\par\medskip
   Assume that $b\leq m_{*}-31$, then $E_{k}\in H^{b+1}(X)$ by
Theorem 2.3.  The following estimate of $E_{k}$ will be utilized
in the next proposition:
\begin{equation}
\parallel E_{k}\parallel_{b+1}\leq\sum_{i=0}^{k-1}\parallel
e_{i}\parallel_{b+1}\leq\varepsilon\delta^{2}\sum_{i=0}^{k-1}
\mu_{i}\leq\varepsilon(\sum_{i=0}^{\infty}\mu^{-1}_{i})\delta^{2}\mu^{k}
\leq\varepsilon(\sum_{i=0}^{\infty}2^{-i})\delta^{2}\mu_{k}.
\end{equation}

\textbf{Proposition 4.3.}  \textit{If the hypotheses of
Proposition} 4.2 \textit{hold and $b\leq m_{*}-31$, then}
VII$_{k+1}$ \textit{holds for all $0\leq m\leq m_{*}$.}\medskip

\textit{Proof.}  By Lemma 3.1 $(iii)$,
\begin{eqnarray}
& &\!\!\!\parallel f_{k+1}\parallel_{m,\text{
}\mathbb{R}^{2}}\\
&\leq&\!\!\!\varepsilon^{-1} \!\parallel T\parallel
C_{m,19}(\parallel
S_{k}E_{k}-S_{k+1}E_{k+1}+(S_{k}-S_{k+1})\Phi(w_{0})\parallel_{m}\nonumber\\
& &\!\!\!+\parallel v_{k+1}\parallel_{m+4}\parallel
S_{k}E_{k}-S_{k+1}E_{k+1}+(S_{k}-S_{k+1})\Phi(w_{0})\parallel_{2}).\nonumber
\end{eqnarray}
Furthermore using (54) and the estimate
$\parallel\Phi(w_{0})\parallel_{b+1}\leq C_{b}\varepsilon^{2n}$,
we obtain for all $m\geq b+1$,
\begin{eqnarray}
& &\!\!\!\parallel
S_{k}E_{k}-S_{k+1}E_{k+1}+(S_{k}-S_{k+1})\Phi(w_{0})\parallel_{m}\\
&\leq& \!\!\!C_{m,20}(\mu_{k}^{m-b-1}\parallel
E_{k}\parallel_{b+1}+\mu_{k+1}^{m-b-1}\parallel
E_{k+1}\parallel_{b+1}\nonumber\\
& &\!\!\!+(\mu_{k}^{m-b-1}+\mu_{k+1}^{m-b-1})\parallel\Phi(w_{0})\parallel_{b+1})\nonumber\\
&\leq&\!\!\!
C_{m,21}\varepsilon\delta^{2}(1+\mu^{b-m})\mu_{k+1}^{m-b}.\nonumber
\end{eqnarray}
If $m<b+1$, then applying similar methods along with VI$_{k+1}$ to
\begin{eqnarray*}
& &\parallel
S_{k}E_{k}-S_{k+1}E_{k+1}+(S_{k}-S_{k+1})\Phi(w_{0})\parallel_{m}\\
&\leq& \parallel (I-S_{k})E_{k}\parallel_{m}+\parallel
(I-S_{k+1})E_{k}\parallel_{m}+\parallel
S_{k+1}e_{k}\parallel_{m}\\
& &+\parallel (I-S_{k})\Phi(w_{0})\parallel_{m}+\parallel
(I-S_{k+1})\Phi(w_{0})\parallel_{m},
\end{eqnarray*}
yields the same estimate found in (56).  Therefore plugging into
(55) produces
\begin{eqnarray*}
\parallel f_{k+1}\parallel_{m,\text{ }\mathbb{R}^{2}}&\leq&
C_{m,22}[\delta^{2}(1+\mu^{b-m})\mu^{m-b}_{k+1}
+\delta^{3}(1+\mu^{b-2})\mu^{m+6-2b}_{k+1}]\\
&\leq& C_{m,23}\delta^{2}(1+\mu^{b-m})\mu^{m-b}_{k+1},
\end{eqnarray*}
if $m+4-b\geq 1/2$.  If $m+4-b\leq-1/2$ and $m\geq 2$, then using
$\parallel v_{k+1}\parallel_{m+4}\leq C_{3}\delta$ in the estimate
above gives the desired result. Moreover if $0\leq m<2$, then in
place of (55) we use the estimate
\begin{equation*}
\parallel f_{k+1}\parallel_{m,\text{ }\mathbb{R}^{2}}\leq\varepsilon^{-1}\parallel
T\parallel C_{m,24}\parallel
S_{k}E_{k}-S_{k+1}E_{k+1}+(S_{k}-S_{k+1})\Phi(w_{0})\parallel_{m}
\end{equation*}
combined with the above method to obtain the desired result.
Lastly if $m+4-b=0$, then replace $\parallel
v_{k+1}\parallel_{m+4}$ in (55) by $\parallel
v_{k+1}\parallel_{m+5}$ and follow the above method.
q.e.d.\medskip

\textbf{Proposition 4.4.}  \textit{If the hypotheses of
Proposition} 4.3 \textit{hold and $b=m_{*}-31$, then} VIII$_{k+1}$
\textit{holds for $0\leq m\leq m_{*}-30$.}\medskip

\textit{Proof.}  By (53), VI$_{k+1}$, and $m\leq b+1=m_{*}-30$, we
have
\begin{eqnarray*}
& &\!\!\!\parallel\Phi(w_{k+1})\parallel_{m}\\
&\leq&\!\!\!\parallel(I-S_{k})\Phi(w_{0})\parallel_{m}+\parallel(I-S_{k})E_{k}\parallel
_{m}+\parallel e_{k}\parallel_{m}\\
&\leq&\!\!\!
C_{m,25}(\mu^{m-b-1}_{k}\parallel\Phi(w_{0})\parallel_{b+1}+
\mu^{m-b-1}_{k}\parallel
E_{k}\parallel_{b+1}+\varepsilon\delta^{2}\mu^{m-b}_{k}).
\end{eqnarray*}
Applying the estimate (54),
$\parallel\Phi(w_{0})\parallel_{b+1}\leq C_{b}\varepsilon^{2n}\leq
\delta^{2}$, and $\delta\mu^{b-m}\leq\varepsilon^{\frac{1}{b+1}}$
produces
\begin{equation*}
\parallel\Phi(w_{k+1})\parallel_{m}\leq
C_{m,26}(\delta^{2}\mu^{b-m}+\varepsilon\delta^{2}
\mu^{b-m})\mu_{k+1}^{m-b}\leq\delta\mu_{k+1}^{m-b},
\end{equation*}
if $\varepsilon$ is sufficiently small.  q.e.d.\par\medskip
   To complete the proof by induction we will now prove the case
$k=0$.  Since $w_{0}=0$, II$_{0}$, III$_{0}$, IV$_{0}$, and
V$_{0}$ are trivial.  Furthermore since
$\parallel\Phi(w_{0})\parallel_{m}\leq\varepsilon\delta^{2}$ if
$\varepsilon=\varepsilon(m)$ is sufficiently small and $m\leq
m_{*}$, VII$_{0}$ and VIII$_{0}$ hold.  In addition, by Lemma 4.2
we can apply Theorem 3.2 to obtain
\begin{equation*}
\parallel u_{0}\parallel_{m}\leq C_{m}\parallel
f_{0}\parallel_{m,\text{ }\mathbb{R}^{2}}\leq
C^{'}_{m}\delta^{2}\leq\delta
\end{equation*}
if $\delta$ is small, so that I$_{1}$ is valid.  Lastly, the proof
of Proposition 4.2 now shows that VI$_{1}$ is valid. This
completes the proof by induction.\par
   In view of the hypotheses of Propositions 4.1-4.4, we require
$m_{*}\geq 58$ and choose $b=m_{*}-31$.  The following corollaries
will complete the proof of Theorem 0.3.\medskip

\textbf{Corollary 4.1.}  $w_{k}\rightarrow w$\textit{ in
}$H^{m_{*}-32}(X)$.\medskip

\textit{Proof.}  For $0\leq m\leq m_{*}-32$ and $i>j$, I$_{k}$
implies that
\begin{equation*}
\parallel w_{i}-w_{j}\parallel_{m}\leq
\sum_{k=j}^{i-1}\parallel u_{k}\parallel_{m}\leq
\delta\sum_{k=j}^{i-1}\mu_{k}^{m-b}\leq\delta
\sum_{k=j}^{i-1}\mu^{-k}.
\end{equation*}
Hence, $\{w_{k}\}$ is Cauchy in $H^{m}(X)$ for all $0\leq m\leq
m_{*}-32$.  q.e.d.\medskip

\textbf{Corollary 4.2.}  $\Phi(w_{k})\rightarrow 0$\textit{ in
}$C^{0}(X)$.\medskip

\textit{Proof.}  By the Sobolev Lemma and VIII$_{k}$,
\begin{equation*}
|\Phi(w_{k})|_{C^{0}(X)}\leq
C\parallel\Phi(w_{k})\parallel_{2}\leq C\delta\mu^{2-b}_{k}.
\end{equation*}
The desired conclusion follows since $b=m_{*}-31\geq 27$.
q.e.d.\par\medskip
   Let $r,K,a_{ij},$ and $f$ be as in Theorem 0.3.  If
$K,a_{ij},f\in C^{r}$, $r\geq 58$, then there exists a $C^{r-34}$
solution of (50).\bigskip

\textbf{Remark.}  After completion of this manuscript, it was
brought to the author's attention that the methods of [4] and [7]
may be adapted to help simplify the linear existence theory of
sections $\S 1$ and $\S 2$.\bigskip

\textbf{Acknowledgments.}  This is a revised portion of my
dissertation [9] conducted at the University of Pennsylvania under
the direction of Professor Jerry Kazdan.  I would like to thank
Jerry Kazdan, Dennis DeTurck, Herman Gluck, and Stephen Shatz for
their suggestions and assistance.  Also a special thanks for very
useful discussions is due to Qing Han, who has obtained a similar
result [5] independently for the isometric embedding problem,
Theorem 0.1.\pagebreak
\bigskip
\begin{center}
\textbf{5.  Appendix}
\end{center}
\bigskip

   Here we shall show that Theorem 0.1 holds for an arbitrary
smooth curve $\sigma$ passing through the origin.  This will be
accomplished by utilizing the special structure of the isometric
embedding equation (1), to show that the calculations of Lemma 1.2
can be refined in this case so that the canonical form (4) may be
achieved without requiring the Christoffel symbols to vanish along
$\sigma$.  This observation is due to Qing Han.  Recall that the
geodesic hypothesis on $\sigma$ was only used to obtain a high
degree of vanishing for the Christoffel symbols along
$\sigma$.\par
   Let $g=g_{ij}du_{i}du_{j}$ be the given metric in local
coordinates, and write equation (1) as
\begin{equation*}
\det\nabla_{ij}z=K|g|(1-|\nabla_{g}z|^{2}),
\end{equation*}
where $\nabla_{ij}$ are covariant derivatives, $K$ is the Gaussian
curvature, $\nabla_{g}$ is the gradient operator with respect to
$g$, and $|g|=\det g_{ij}$.  Following the set up of the
introduction we set $u_{i}=\varepsilon^{2}x_{i}$, and
$z=u_{1}^{2}/2+\varepsilon^{5}w$.  Then as in (7) the
linearization of (6) becomes
\begin{equation*}
L_{1}(w)v=\sum_{i,j}b^{ij}v_{;ij}+\sum_{i}b^{i}v_{;i}:=
\sum_{i,j}b_{ij}^{1}v_{x_{i}x_{j}}+\sum_{i}b_{i}^{1}v_{x_{i}},
\end{equation*}
where $v_{;ij}$, $v_{;i}$ denote covariant derivatives in $x_{i}$
coordinates (we will denote covariant derivatives in $u_{i}$
coordinates by $\nabla_{ij}v$, $\nabla_{i}v$), $b^{ij}$ is the
cofactor matrix given by
\begin{eqnarray}
b^{11}=b_{11}^{1}\!\!&=&\!\!\varepsilon\nabla_{22}z=
\varepsilon^{2}O(1+|\nabla w|+|\nabla^{2}w|),\\
b^{12}=b^{21}=b_{12}^{1}=b_{21}^{1}\!\!&=&\!\!
-\varepsilon\nabla_{12}z=\varepsilon^{2}O(1+|\nabla w|+|\nabla^{2}w|),\nonumber\\
b^{22}=b_{22}^{1}\!\!&=&\!\!\varepsilon\nabla_{11}z=\varepsilon(1+\varepsilon
O(1+|\nabla w|+|\nabla^{2}w|)), \nonumber
\end{eqnarray}
and
\begin{equation*}
b_{i}^{1}=-b^{lk}\Gamma_{lk}^{i}+b^{i}=-b^{lk}\Gamma_{lk}^{i}+
\varepsilon^{2(n+1)}H^{n+1}(x_{1},x_{2})P_{i}(\varepsilon,x_{1},x_{2},
\nabla w)
\end{equation*}
for some $P_{i}$, with $\Gamma_{lk}^{i}$ Christoffel symbols for
$g$ in $x_{i}$ coordinates.  Also throughout this section the
summation convention for raised and lowered indices will be
used.\par
   Continuing to follow the procedure of section $\S 1$, we find
that (9) produces
\begin{equation*}
L_{3}(w)v=\sum_{i,j}b_{ij}^{3}v_{x_{i}x_{j}}+\sum_{i}b_{i}^{3}v_{x_{i}}
\end{equation*}
where
\begin{eqnarray}
b_{11}^{3}\!\!&=&\!\!(b^{22})^{-2}((b^{12})^{2}+\varepsilon^{2(n+2)}H^{n+1}P),\nonumber\\
b_{12}^{3}\!\!&=&\!\!(b^{22})^{-1}b^{12},\\
b_{22}^{3}\!\!&=&\!\!1,\nonumber\\
b_{i}^{3}\!\!&=&\!\!(b^{22})^{-1}(-b^{lk}\Gamma_{lk}^{i}
+\varepsilon^{2(n+1)}H^{n+1}P_{i}),\nonumber
\end{eqnarray}
and
\begin{equation*}
K|g|(1-|\nabla_{g}z|^{2})=\varepsilon^{2(n+1)}
H^{n+1}(x_{1},x_{2})P(\varepsilon,x_{1},x_{2},\nabla w).
\end{equation*}
Let
\begin{equation*}
\xi=\xi(x_{1},x_{2}),\text{ }\text{ }\text{ }\text{ }\text{
}\eta=x_{2},
\end{equation*}
be the change of coordinates of Lemma 1.2, so that $\xi$ satisfies
(10):
\begin{equation}
b^{12}\xi_{x_{1}}+b^{22}\xi_{x_{2}}=0.
\end{equation}
If as before $b_{ij}^{4}$ and $b_{i}^{4}$ denote the coefficients
of $L_{3}(w)$ in these new coordinates, then all the conclusions
of Lemma 1.2 hold.  In fact the proof requires no modification,
except to justify the expression for $b_{1}^{4}$ which we now
show.\par
   Using (13), (14), and (58) we obtain the analogue of (15):
\begin{eqnarray}
b_{1}^{4}\!\!&=&\!\!\sum_{i,j}b_{ij}^{3}\xi_{x_{i}x_{j}}+\sum_{i}
b_{i}^{3}\xi_{x_{i}}\nonumber\\
&=&\!\!\frac{\varepsilon^{2(n+2)}H^{n+1}P}{(b^{22})^{2}}\xi_{x_{1}x_{1}}
-\left[\left(\frac{b^{12}}{b^{22}}\right)\left(\frac{b^{12}}{b^{22}}\right)_{x_{1}}
+\left(\frac{b^{12}}{b^{22}}\right)_{x_{2}}\right]\xi_{x_{1}}\\
& &\!\!+\sum_{i}b_{i}^{3} \xi_{x_{i}}.\nonumber
\end{eqnarray}
Calculating the second term on the right-hand side of (60) yields,
\begin{eqnarray*}
&
&\!\!(b^{22})^{2}\left[\left(\frac{b^{12}}{b^{22}}\right)\left(\frac{b^{12}}{b^{22}}\right)_{x_{1}}
+\left(\frac{b^{12}}{b^{22}}\right)_{x_{2}}\right]\\
\!\!&=&\!\!
b^{12}b^{12}_{x_{1}}-(b^{22})^{-1}(b^{12})^{2}b^{22}_{x_{1}}+b^{22}b^{12}_{x_{2}}
-b^{12}b^{22}_{x_{2}}\\
&=&\!\!
b^{12}b^{12}_{x_{1}}-b^{11}b^{22}_{x_{1}}+b^{22}b^{12}_{x_{2}}-b^{12}b^{22}_{x_{2}}\\
& &\!\!
+(b^{22})^{-1}b^{22}_{x_{1}}(\det b^{ij})\\
&=&\!\!-b^{12}_{x_{1}}b^{12}+b^{11}_{x_{1}}b^{22}+b^{22}b^{12}_{x_{2}}
-b^{12}b^{22}_{x_{2}}\\
& &\!\!+(b^{22})^{-1}b^{22}_{x_{1}}(\det b^{ij})-(\det
b^{ij})_{x_{1}}.
\end{eqnarray*}
Therefore, (59) and (60) imply that
\begin{eqnarray}
b^{22}b_{1}^{4}\!\!&=&\!\!-(b_{x_{1}}^{11}+b_{x_{2}}^{12}+b^{lk}\Gamma_{lk}^{1}
-\varepsilon^{2(n+1)}H^{n+1}P_{1}-((b^{22})^{-1}\det b^{ij})_{x_{1}})\xi_{x_{1}}\nonumber\\
&
&\!\!-(b_{x_{1}}^{12}+b^{22}_{x_{2}}+b^{lk}\Gamma_{lk}^{2}-\varepsilon^{2(n+1)}
H^{n+1}P_{2})\xi_{x_{2}}\\
&
&\!\!+(b^{22})^{-1}\varepsilon^{2(n+2)}H^{n+1}P\xi_{x_{1}x_{1}}.\nonumber
\end{eqnarray}
Lastly, from (57) we calculate
\begin{eqnarray}
&
&\!\!\!\varepsilon^{3}(b^{11}_{x_{1}}+b^{12}_{x_{2}}+b^{lk}\Gamma_{lk}^{1})\\
&=&\!\!
-\Gamma_{j2}^{j}z_{x_{1}x_{2}}+\Gamma_{j1}^{j}z_{x_{2}x_{2}}\nonumber\\
&+&\!\!\!(\Gamma_{12,x_{2}}^{i}
-\Gamma_{22,x_{1}}^{i}-\Gamma_{11}^{1}\Gamma_{22}^{i}+2\Gamma_{12}^{1}\Gamma_{12}^{i}
-\Gamma_{22}^{1}\Gamma_{11}^{i})z_{x_{i}}\nonumber\\
&=&\!\!\!\varepsilon^{3}\Gamma_{j2}^{j}b^{12}+\varepsilon^{3}\Gamma_{j1}^{j}b^{11}\nonumber\\
&+&\!\!\!(\Gamma_{12,x_{2}}^{i}\!
-\Gamma_{22,x_{1}}^{i}\!-\Gamma_{11}^{1}\Gamma_{22}^{i}\!+2\Gamma_{12}^{1}\Gamma_{12}^{i}\!
-\Gamma_{22}^{1}\Gamma_{11}^{i}\!-\Gamma_{j2}^{j}\Gamma_{12}^{i}\!
+\Gamma_{j1}^{j}\Gamma_{22}^{i})z_{x_{i}}.\nonumber
\end{eqnarray}
However, we see that the coefficient of $z_{x_{i}}$ is in fact a
curvature term.  More precisely, if we denote it by $\Omega^{i}$
then
\begin{eqnarray}
\Omega^{i}\!\!&=&\!\!\Gamma_{12,x_{2}}^{i}-\Gamma_{22,x_{1}}^{i}
+\Gamma_{12}^{j}\Gamma_{j2}^{i}-\Gamma_{22}^{j}\Gamma_{j1}^{i}\\
&=&\!\!-\varepsilon^{4}R^{i}_{212}=-\varepsilon^{4}g^{i1}|g|K=-\varepsilon^{2(n+3)}
H^{n+1}\overline{P}_{1}^{i}\nonumber
\end{eqnarray}
for some $\overline{P}_{1}^{i}$, where $R^{i}_{jkl}$ is the
Riemann tensor for $g$ in $u_{i}$ coordinates (recall that
$\Gamma_{lk}^{i}$ are Christoffel symbols in $x_{i}$ coordinates).
A similar calculation shows that
\begin{equation}
\varepsilon^{3}(b^{12}_{x_{1}}+b^{22}_{x_{2}}+b^{lk}\Gamma_{lk}^{2})
=\varepsilon^{3}\Gamma_{j1}^{j}b^{12}+\varepsilon^{3}\Gamma_{j2}^{j}b^{22}
-\varepsilon^{2(n+3)}H^{n+1}\overline{P}_{2}^{i}z_{x_{i}}
\end{equation}
for some $\overline{P}_{2}^{i}$.  Then observing that
\begin{equation}
\det b^{ij}=\varepsilon^{2}\Phi(w)+\varepsilon^{2(n+2)}H^{n+1}P
\end{equation}
from (6), we may combine (59) and (61)-(65) to obtain the desired
expression for $b_{1}^{4}$ as stated in Lemma 1.2 (note that the
linear combination of $\Phi(w)$ and $\partial_{x_{1}}\Phi(w)$ will
appear slightly different than in Lemma 1.2).  Having established
Lemma 1.2, we can then apply the remainder of section $\S 1$ as
well as sections $\S 2$, $\S 3$, and $\S 4$ without change in
order to obtain Theorem 0.1 for an arbitrary smooth curve
$\sigma$.
\bigskip
\begin{center}
\textbf{References}
\end{center}
\bigskip

\noindent[1]\hspace{.06in} G. Birkhoff, G.-C. Rota,
\textit{Ordinary Differential Equations}, Blaisdell
\par\hspace{.06in} Publishing, London, 1969, MR 0972977, Zbl 0377.34001.\bigskip

\noindent[2]\hspace{.06in} K. O. Friedrichs, \textit{The identity
of weak and strong extensions of}
\par\hspace{.06in} \textit{differential operators,} Trans. Amer.
Math. Soc., \textbf{55} (1944), 132- \par\hspace{.06in} 151, MR
0009701, Zbl 0061.26201.
\bigskip

\noindent[3]\hspace{.06in} S. Gallerstedt, \textit{Quelques
probl\`{e}mes mixtes pour l'\'{e}quation} $y^{m}z_{xx}+$
\par\hspace{.06in} $z_{yy}=0$, Arkiv f\"{o}r Matematik, Astronomi och Fysik,
\textbf{26A} (1937), \par\hspace{.06in} no. 3, 1-32.\bigskip

\noindent[4]\hspace{.06in} Q. Han, \textit{On the isometric
embedding of surfaces with Gauss curva-} \par\hspace{.06in}
\textit{ture changing sign cleanly}, Comm. Pure Appl. Math.,
\textbf{58} (2005), \par\hspace{.06in} 285-295, MR 2094852, Zbl
1073.53005.\bigskip

\noindent[5]\hspace{.06in} Q. Han, \textit{Local isometric
embedding of surfaces with Gauss curvature}
\par\hspace{.06in} \textit{changing sign stably across a curve},
Cal. Var. \& P.D.E., \textbf{25} (2006), \par\hspace{.06in} no. 1,
79--103, MR 2183856, Zbl pre05009621.\bigskip

\noindent[6]\hspace{.06in} Q. Han, J.-X. Hong, C.-S. Lin,
\textit{Local isometric embedding of sur-} \par\hspace{.06in}
\textit{faces with nonpositive Gaussian curvature,} J.
Differential Geom.,
\par\hspace{.06in} \textbf{63} (2003), 475-520, MR 2015470, Zbl 1070.53034.\bigskip

\noindent[7]\hspace{.06in} J.-X. Hong, \textit{Cauchy problem for
degenerate hyperbolic Monge-Amp\`{e}re}
\par\hspace{.06in} \textit{equations}, J. Partial Diff. Equations,
\textbf{4}
(1991), 1-18, MR 1111376.\bigskip

\noindent[8]\hspace{.06in} H. Jacobowitz, \textit{Local isometric
embeddings}, Seminar on Differential \par\hspace{.06in} Geometry,
Annals of Math. Studies, \textbf{102}, edited by S.-T. Yau, 1982,
\par\hspace{.06in} 381-393, MR 0645749, Zbl 0481.53018.\bigskip

\noindent[9]\hspace{.09in} M. A. Khuri, \textit{The local
isometric embedding in $\mathbb{R}^{3}\!$ \hspace{-.02in} of
two-dimensional}
\par\hspace{.06in} \textit{Riemannian manifolds with Gaussian
curvature changing sign to} \par\hspace{.06in} \textit{finite
order on a curve}, Dissertation, University of Pennsylvania,
\par\hspace{.06in} 2003.\bigskip

\noindent[10]  M. A. Khuri, \textit{Counterexamples to the local
solvability of Monge-} \par\hspace{.06in} \textit{Amp\`{e}re
equations in the plane,} Comm. PDE, $\mathbf{32}$ (2007),
665-674.\medskip

\noindent[11]  M. A. Khuri, \textit{Local solvability of
degenerate Monge-Amp\`{e}re equa-}
\par\hspace{.06in} \textit{tions and applications to
geometry}, Electron. J. Diff. Eqns., $\mathbf{2007}$
\par\hspace{.06in} (2007), No. 65, 1-37.\bigskip

\noindent[12]  C.-S. Lin, \textit{The local isometric embedding in
$\mathbb{R}^{3}$ of} 2-\textit{dimensional} \par\hspace{.06in}
\textit{Riemannian manifolds with nonnegative curvature}, J.
Differential
\par\hspace{.06in} Geom., \textbf{21} (1985), no. 2, 213-230, MR 0816670,
Zbl 0584.53002.\bigskip

\noindent[13]  C.-S. Lin, \textit{The local isometric embedding in
$\mathbb{R}^{3}$ of two-dimensional}
\par\hspace{.06in} \textit{Riemannian manifolds with Gaussian
curvature changing sign cleanly}, \par\hspace{.06in} Comm. Pure
Appl. Math., \textbf{39} (1986), no. 6, 867-887, MR 0859276,
\par\hspace{.06in} Zbl 0612.53013.\bigskip

\noindent[14]  N. Nadirashvili, Y. Yuan, \textit{Improving
Pogorelov's isometric embed-} \par\hspace{.06in} \textit{ding
counterexample}, preprint.\bigskip

\noindent[15]  A. V. Pogorelov, \textit{An example of a
two-dimensional Riemannian} \par\hspace{.06in} \textit{metric not
admitting a local realization in $E_{3}$}, Dokl. Akad. Nauk.
\par\hspace{.06in} USSR, \textbf{198} (1971), 42-43, MR 0286034, Zbl 0232.53013.\bigskip

\noindent[16]  E. G. Poznyak, \textit{Regular realization in the
large of two-dimensional} \par\hspace{.06in} \textit{metrics of
negative curvature}, Soviet Math. Dokl., \textbf{7} (1966), 1288-
\par\hspace{.06in} 1291, MR 0205204, Zbl 0168.19501.\bigskip

\noindent[17]  E. G. Poznyak, \textit{Isometric immersions of
two-dimensional Rie-}
\par\hspace{.06in} \textit{mannian metrics in Euclidean
space}, Russian Math. Surveys, \textbf{28} \par\hspace{.06in}
(1973), 47-77.\bigskip

\noindent[18]  J. T. Schwartz, \textit{Nonlinear Functional
Analysis}, New York Univer- \par\hspace{.06in} sity, New York,
1964, MR 0433481, Zbl 0203.14501.\bigskip

\noindent[19]  E. Stein, \textit{Singular Integrals and
Differentiability Properties of Func-}
\par\hspace{.06in} \textit{tions}, Princeton University Press,
Princeton, 1970, MR 0290095, \par\hspace{.06in} Zbl
0207.13501.\bigskip

\noindent[20]  M. E. Taylor, \textit{Partial Differential
Equations} III, Springer-Verlag, \par\hspace{.06in} New York,
1996,  MR 1477408, Zbl 0869.35004.\bigskip

\noindent[21]  J. Weingarten, \textit{\"{U}ber die theorie der
Aubeinander abwickelbarren} \par\hspace{.06in}
\textit{Oberfl\"{a}chen}, Berlin, 1884.\par

\bigskip\bigskip\bigskip
\begin{small}
\noindent\textsc{Stanford University}\\
\hspace{4.0in}\textsc{Mathematics, Bldg. 380}\\
\hspace{4.0in}\textsc{450 Serra Mall}\\
\hspace{4.0in}\textsc{Stanford, CA 94305-2125}\\
\hspace{4.0in}\textit{E-mail address:} khuri@math.stanford.edu
\end{small}

\end{document}